\newif\ifpdf
\ifx\pdfoutput\undefined
\pdffalse % we are not running PDFLaTeX
\else
\pdftrue \fi

\newif\iffinal
\finaltrue % Final version

\documentclass[reqno,twoside,11pt]{amsart}
\iffinal\else\usepackage[notref,notcite]{showkeys}\fi

\usepackage{amsmath}
\usepackage{amsfonts}
\usepackage{amssymb}
\usepackage{color}
\usepackage{multido}
\usepackage{pstricks,pst-plot,pstricks-add,pst-math}
\usepackage{enumerate}
\usepackage{verbatim}
\usepackage{epsfig}
\usepackage{setspace}

\IfFileExists{epsf.def}{\input epsf.def}{\usepackage{epsf}}

\IfFileExists{myowntimes.sty}{\usepackage{myowntimes}\usepackage{mathrsfs}}
        {\usepackage{mathpazo}\usepackage{mathrsfs}}

\DeclareFontFamily{OT1}{eusb}{} \DeclareFontShape{OT1}{eusb}{m}{n}
{<5> <6> <7> <8> <9> <10> <11> <12> <14.4> eusb10}{}
\DeclareMathAlphabet{\eusb}{OT1}{eusb}{m}{n}

\DeclareFontFamily{OT1}{eusm}{} \DeclareFontShape{OT1}{eusm}{m}{n}
{<5> <6> <7> <8> <9> <10> <11> <12> <14.4> eusm10}{}
\DeclareMathAlphabet{\eusm}{OT1}{eusm}{m}{n}

\DeclareFontFamily{OT1}{eufm}{} \DeclareFontShape{OT1}{eufm}{m}{n}
{<5> <6> <7> <8> <9> <10> <11> <12> <14.4> eufm10}{}
\DeclareMathAlphabet{\mathfrak}{OT1}{eufm}{m}{n}

\DeclareFontFamily{OT1}{fraktura}{}
\DeclareFontShape{OT1}{fraktura}{m}{n} {<5> <6> <7> <8> <9> <10> <11>
  <12> <13> <14.4> [1.1] eufm10}{}
\DeclareMathAlphabet{\fraktura}{OT1}{fraktura}{m}{n}

\DeclareFontFamily{OT1}{cmfi}{} \DeclareFontShape{OT1}{cmfi}{m}{n}
{<5> <6> <7> <8> <9> <10> <11> <12> <13> <14.4> [0.9] cmfi10}{}
\DeclareMathAlphabet{\cmfi}{OT1}{cmfi}{b}{n}

\DeclareFontFamily{OT1}{cmss}{} \DeclareFontShape{OT1}{cmss}{m}{n}
{<5> <6> <7> <8> <9> <10> <11> <12> <13> <14.4> cmss10}{}
\DeclareMathAlphabet{\cmss}{OT1}{cmss}{m}{n}
\newtheoremstyle{thm}{1.5ex}{1.5ex}{\itshape\rmfamily}{}
{\bfseries\rmfamily}{}{2ex}{}

\newtheoremstyle{def}{1.5ex}{1.5ex}{\slshape\rmfamily}{}
{\bfseries\rmfamily}{}{2ex}{}

\newtheoremstyle{rem}{1.3ex}{1.3ex}{\rmfamily}{}
{\itshape}
{} {1.5ex}{}

\theoremstyle{thm}
\newtheorem{maintheorem}{Theorem}

\newtheorem{theorem}{Theorem}[section]
\newtheorem{lemma}[theorem]{Lemma}

\newtheorem*{Main Theorem}{Main Theorem.}
\newtheorem{corollary}[theorem]{Corollary}
\newtheorem*{special theorem}{Lindeberg-Feller Theorem for Martingales}

\newtheorem{definition}[theorem]{Definition}

\theoremstyle{rem}

\numberwithin{equation}{section}

\renewcommand{\section}{\secdef\sct\sect}
\newcommand{\sct}[2][default]{%
\refstepcounter{section}
\addcontentsline{toc}{section}{{\tocsection
    {}{\thesection}{\!\!\!\!#1\dotfill}}{}}
\vspace{0.7cm}
\centerline{\scshape\thesection.\ #1} \nopagebreak \vspace{0.2cm}}
\newcommand{\sect}[1]{%
\vspace{0.4cm} \centerline{\large\scshape\rmfamily #1}
\vspace{0.2cm}
}

\renewcommand{\subsection}{\secdef\subsct\sbsect}
\newcommand{\subsct}[2][default]{\refstepcounter{subsection}
\addcontentsline{toc}{subsection}
{{\tocsection{\!\!}{\hspace{1.2em}\thesubsection}{\!\!\!\!#1\dotfill}}{}}
\nopagebreak
{\flushleft\bf
\thesubsection~\bf #1.~}
\noindent
\nopagebreak}
\newcommand{\sbsect}[1]{
\noindent
\textbf{#1.~}
}

\renewcommand{\subsubsection}{%
\secdef \subsubsect\sbsbsect}
\newcommand{\subsubsect}[2][default]{%
\refstepcounter{subsubsection}
\addcontentsline{toc}{subsubsection}{{\tocsection{\!\!}
{\hspace{3.05em}\thesubsubsection}{\!\!\!\!#1\dotfill}}{}}
\nopagebreak
\vspace{0.15\baselineskip} \nopagebreak {\flushleft\rmfamily
\itshape\thesubsubsection
\ \rmfamily #1\/.}\ }
\newcommand{\sbsbsect}[1]{\vspace{0.1cm}\noindent
\rmfamily \itshape
\arabic{section}.\arabic{subsection}.\arabic{subsubsection} \
\sffamily #1\/.\ }

\iffinal

\else

\fi

\renewcommand{\caption}[1]{%
\vglue0.5cm
\refstepcounter{figure}
\begin{minipage}{0.9\textwidth}\small {\sc Figure~\thefigure. }#1\end{minipage}}

\newcommand{\ee}{\end{equation}}
\newcommand{\be}{\begin{equation}}
\newcommand{\eml}{\end{multline}}
\newcommand{\bml}{\begin{multline}}
\newcommand{\ra}{\rightarrow}

 \newcommand{\set}[1]{\left\{#1\right\}}

\def\qed{ \hfill $\square$}

 \DeclareMathOperator{\E}{\mathbb{E}}

 \newcommand{\bN}{\mathbb{N}}

 \newcommand{\mvone}{\boldsymbol{1}}

\newcommand{\textd}{\text{\rm d}\mkern0.5mu}

\renewcommand{\AA}{\mathcal A}
\newcommand{\BB}{\mathcal B}
\newcommand{\CC}{\mathcal C}
\newcommand{\DD}{\mathcal D}
\newcommand{\EE}{\mathcal E}
\newcommand{\FF}{\mathcal F}
\newcommand{\GG}{\mathcal G}
\newcommand{\HH}{\mathcal H}
\newcommand{\II}{\mathcal I}

\newcommand{\LL}{\mathcal L}

\newcommand{\NN}{\mathcal N}
\newcommand{\OO}{\mathcal O}
\newcommand{\PP}{\mathcal P}

\newcommand{\RR}{\mathcal R}

\newcommand{\ZZ}{\mathcal Z}

\newcommand{\BbbP}{\mathbb P}

\newcommand{\R}{\mathbb R}

\newcommand{\Z}{\mathbb Z}

\newcommand{\scrC}{\mathscr{C}}

\newcommand{\scrE}{\mathscr{E}}

\newcommand{\scrG}{\mathscr{G}}

\newcommand{\scrM}{\mathscr{M}}

\newcommand{\scrO}{\mathscr{O}}

\newcommand{\scrR}{\mathscr{R}}
\newcommand{\scrS}{\mathscr{S}}

\newcommand{\scrX}{\mathscr{X}}

\def\myffrac#1#2 in #3{\raise 2.6pt\hbox{$#3 #1$}\mkern-1.5mu\raise 0.8pt\hbox{$#3/$}\mkern-1.1mu\lower 1.5pt\hbox{$#3 #2$}}
\newcommand{\ffrac}[2]{\mathchoice%
        {\myffrac{#1}{#2} in \scriptstyle}
        {\myffrac{#1}{#2} in \scriptstyle}
        {\myffrac{#1}{#2} in \scriptscriptstyle}
        {\myffrac{#1}{#2} in \scriptscriptstyle}
}

\newcommand{\la}{\langle}
\newcommand{\rra}{\rangle}
\title{Heat Kernel Upper Bounds on Long Range Percolation Clusters}
\author[N. Crawford]{Nicholas Crawford$^1$}
\author[A. Sly]{Allan Sly$^2$}
\thanks{{\tt
    email:allansly@gmail.com}}  
\thanks{{\tt
    email:nickcrawford12345@gmail.com}, Supported in part by DOD ONR grant
  N0014-07-1-05-06}
\thanks{Work completed while the authors were at UC Berkeley Department of Statistics}
\begin{document}
\thanks{\hglue-4.5mm\fontsize{9.6}{9.6}\selectfont\copyright\,2009 by N.~Crawford and A.~Sly. Reproduction, by any means, of the entire
article for non-commercial purposes is permitted without charge.\vspace{2mm}}
\maketitle
\centerline{$^1$  \textit{Department of Industrial Engineering, The Technion; Haifa, Israel}}
\centerline{$^2$  \textit{The Theory Group, Microsoft Research; Redmond, Washington}}
\begin{abstract}
In this paper, we derive upper bounds for the heat kernel of the simple random walk on the infinite cluster of a supercritical long range percolation process.  For any $d \geq 1$ and for any exponent $s \in (d, (d+2) \wedge 2d)$ giving the rate of decay of the percolation process, we show that the return probability decays like $t^{-\ffrac{d}{s-d}}$ up to logarithmic corrections, where $t$ denotes the time the walk is run.  Our methods also yield generalized bounds on the spectral gap of the dynamics and on the diameter of the largest component in a box.
 Besides its intrinsic interest, the main result is needed for a companion paper studying the scaling limit of simple random walk on the infinite cluster.  \end{abstract}

\section{Introduction}
Random Walks in Random Environments (RWRE) are a major topic of research in probability theory, with important connections to topics such as Anderson Localization and Interface measures from mathematical physics and the study of random matrices through the RWRE transition kernel.  One simple and fundamental setting for RWRE occurs from random bond dilution of an ambient graph $G$ via an $i.i.d.$ percolation process.  Among the basic questions is: to what extent  such a process affects the macroscopic characteristics of the associated Simple Random Walk (SRW)?

In the context of nearest neighbor percolation on $\Z^d$, many authors have contributed to the understanding of the behavior of the SRW.
{ Let us mention some of the key papers in the area: for various models of RWRE, Kipnis-Varadhan \cite{Kip-Var} introduced "the environment viewed from the particle" point of view to derive annealed functional central limit theorems; this work was strengthened in Demasi et al. \cite{Demasi} where it was applied to SRW on nearest neighbor percolation clusters; Sidoravicious-Sznitman \cite{Sid-Szn} extended the percolation theorem of \cite{Demasi} to the quenched regime on $\Z^d, d \geq 4$; Remy-Mathieu \cite{RM}, Barlow \cite{Barlow-HKP} provide quenched heat kernel bounds on super critical percolation clusters, with earlier estimates obtained by Heicklen and Hoffman \cite{HH}; and finally Mathieu-Pianitskii \cite{Mathieu-Piat}, Berger-Biskup\cite{Berger-Biskup} extend \cite{Sid-Szn} to supercritical clusters on $\Z^d, d\geq 2$.}

In this paper, we consider a variant of this latter problem -- SRW on super critical Long Range Percolation clusters on $\Z^d$ (LRP).  LRP, at least in one dimension, was first considered by Schulman in \cite{Schulman83} and Zhang et. al. \cite{Z}.  It is a random graph process on $\Z^d$ where, independently for each pair of vertices $x, y \in \Z^d$, we attach an edge $\la x, y\rra$ with probability $p_{x, y}$.  Throughout we shall assume (though this is not strictly necessary) an isotropic model:  $p_{x, y} = \texttt{P}(\|x-y\|_2)$ where
\be
\label{Eq;Asymp}
-\frac{\log \texttt{P}(r)}{\log r} \rightarrow s \text{ as $r \ra \infty$}
\ee
for some $s \in \R^+$.  There are a number of transitions for the behavior of this process as a function of $s$, one of which is addressed in this paper. The earliest results of this kind were given by Schulman \cite{Schulman83} and sharpened to the critical case $d=1, s=2$ by Newman-Schulman in \cite{New-Schul}, see also \cite{AKN}.  

More recently, the long range model gained interest in the context of "small world phenomena", see works such as  \cite{Milgram}, \cite{WattsStrog} and \cite{BiskupLR} for discussions.  Benjamini and Berger  \cite{BBDiam} initiated a study of geometric properties of a compactified version of this model, focusing on the asymptotics of the diameter on the discrete cycle $\Z/N\Z$.  Their motivation regarded connections to modeling the topology of the internet, see also \cite{Kleinberg} for a different perspective.  
In the compact setting, \cite{BBDiam} was followed up in \cite{CGS} for higher dimensional discrete tori and more recently by Benjamini, Berger and Yadin \cite{BBY}.  In the last paper, the authors study the mixing time $\tau$ of SRW on $\Z/ N\Z$, providing bounds of the form  
\be
c N^{s-1} \leq \tau \leq C N^{s-1} \log^{\delta}N 
\ee
asymptotically almost surely where $\tau$ is the relaxation time of the SRW.  We will make use of the method of proof from the revised version of the paper \cite{BBYr} which fixes a gap in the published version \cite{BBY}.
We note that all of the above mentioned results for SRW on LRP assume that nearest neighbor connections exist with probability $1$.  

On $\Z^d$, we shall say that the family of probabilities $p_{x, y}=P(\|x-y\|_2)$ is \textit{percolating} if the process admits an infinite connected component.  By translation invariance, this property is a $0-1$ event.  The Burton-Keane argument \cite{Burton-Keane} implies that if $P(\|x-y\|_2)$ percolates, then the infinite component is unique $a.s.$.  It is worth noting however that in our setting, this was already shown in \cite{AKN}.  

Since the work of \cite{Schulman83, New-Schul}, most studies of LRP on $\Z^d$ have focused on conditions for and global aspects of the geometry of percolating clusters.  Berger \cite{BLRP} studied the properties of transience, recurrence of SRW on LRP.  Biskup \cite{BiskupLR, BiskupLR2} studied the graph (or chemical) distance and \textit{diameter} of the super-critical process in the regime $d<s<2d$ providing sharp bounds, up to lower order.

\subsection{Statement of Main Result}

{ Assume that $s \in (d, d+2)$}.
Let us begin by noting that the one step annealed transition probabilities $P_{1}(0, x)$ are proportional to $p_{0,x}$ and thus are heavy tailed and in the domain of attraction of an $\alpha=s-d$ stable law. 
{Given this, a natural assumption is } that the quenched SRW process will asymptotically behave like an $\alpha$-stable L\'{e}vy motion.  

The heat kernel $P_t^\omega{(x,y)}$ should thus be expected to decay {point-wisely} as $ t^{-d/(s-d)}$.  The main result of our paper is to establish a heat kernel upper bound of this order up to logarithmic factors.  Note, however, that this intuition fails in the case of $d=1$ and $s\in(2,3)$ where \cite{BBY} showed that the spectral gap of the dynamics on the torus $\Z/(N\Z)$ is $O(N^2)$ and so behaves instead like Brownian motion.  In a subsequent paper \cite{CS} we will show that the limiting distribution in this case is Brownian motion.

{Let $\EE_{\ZZ_d}$ denote that set of unordered pairs of vertices $\la x, y \rra$ such that $x, y \in \Z^d$.  Let $\Omega= \{0, 1\}^{\EE_{\ZZ^d}}$ denote the sample space for LRP on $Z^d$.  For each $\omega \in \Omega, y \in \Z^d$, let $\textrm{deg}^{\omega}(y)$ denote the degree of the vertex $y$ in $\omega$.}  Let $\mu$ denote the product measure on $\Omega$ defined by the probabilities $(p_{x, y})_{\la x, y\rra \in \EE_{\Z^d}}$.

\begin{maintheorem}
\label{T:HKLRP}
Let us consider $s \in (d, d+2)$ for $d\geq 2$ and $s\in (1, 2)$ for $d=1$.  Assume, for simplicity, that there exists $L$ such that
\be
\label{Eq:Cond}
p_{x,y} = 1- e^{-\beta \|x-y\|^{-s}} \text{ for $\|x-y\|_2 \geq L$}
\ee
and suppose that the $p_{x, y}$ are translation invariant and percolating. Then there exists an event $\Omega_1 \subset \Omega$ with $\mu(\Omega_1)=1$, universal constants $C_1, \delta>0$ and a family of random variables $T_x(\omega)$ with the property that $T_x(\omega)< \infty$ whenever $x,y \in \CC^{\infty}(\omega)$, such that the following holds:
\be
P^{\omega}_{t}(x, y) \leq C_1 \textrm{deg}^{\omega}(y) t^{-d/(s-d)} \log^{\delta} t.
\ee
{for $t \geq T_x(\omega) \vee T_y(\omega)$.}
Moreover on the event $E_x:=\{T_x(\omega)< \infty\}$ and for any $\eta> 0$, there exists $C(\eta)>0$ so that we have
\be
\mu(T_x> k|E_x) \leq  C(\eta)k^{-\eta}
\ee
\end{maintheorem}

This follows the work of Barlow \cite{Barlow-HKP} to the context of long range percolation.  We mention that Kumagai and Misumi \cite{2008KM} derived heat kernel estimates for the LRP graph with  $d=1 , \:s>2$; they behave like the nearest neighbor model.  Note that in this case, it is essential to assume nearest neighbor connections to have an infinite component. 

It is clear from our proof that the same result holds for \textit{any} percolating translation invariant $p(x, y)$ which satisfies \eqref{Eq;Asymp} for $s\in (d, d+2)$, but for convenience of exposition we consider the specific form here.  We emphasize that we do not assume that nearest neighbor edges are occupied with probability one.  This leads to significant technical challenges and much of our work goes into controlling complications which arise in this setting.  For instance we must deal with the fact that some portions of the giant component are only connected via very long edges.  Our results also provide a refinement on the understanding of the geometry of components of the LRP on $\Z^d$, one consequence of which is the following:
\begin{maintheorem}
\label{t:diameter}
With the same conditions as in Theorem \ref{T:HKLRP}, let $\mathcal{C}^1(N)$ denote the largest component of the graph induced by restricting the process to the vertices $B_N(0):=[-N,N]^d$ and let $\mathcal{C}^2(N)$ denote the second largest component.  Then there exists a constant $\delta_2>0$ so that
\be
|\mathcal{C}^2(N)| \leq \log^{\delta_2} N
\ee
asymptotically almost surely.

Further, there exists a constant $\delta$ such that asymptotically almost surely
\[
\hbox{Diameter}(\mathcal{C}^1(N)) \leq \log^\delta N,
\]
where $Diameter$ refers to the graph diameter of the LRP cluster.
Finally, there is $\epsilon>0$ so that
\be
 \mu \left( 0 \leftrightarrow  B^c_{N}(0)| 0\notin \CC^1(N) \right)
 \leq  C N^{- \epsilon}
\ee
where the event $\{0 \leftrightarrow  B^c_{N}(0)\}$ denotes that the cluster containing the origin intersects the complement of $[-N, N]^d$.
\end{maintheorem}
{These statements are proved, respectively, as Lemmas \ref{L:Sizes}, \ref{L:P2} and Corollary \ref{C:connect} This Theorem extends results of {\cite{BBDiam, CGS, BiskupLR2}} to the setting where nearest neighbor bonds do not occur with probability $1$.  It should be noted that \cite{BiskupLR2} obtains the correct power $\delta= \log_2^{-1}(2d/s)$ up to $o(1)$ terms vanishing with $N$ and that all three papers consider $s \in (d, 2d)$, whereas we only explicitly consider $s \in(d, d+2 \wedge 2d)$ (for us this is not really a necessary restriction for \ref{t:diameter} but presumably is for \ref{T:HKLRP} ).  Biskup \cite{BiskupLR} analyses the case without nearest neighbour bonds.  He bounds the distance of two typical vertices but does not bound the diameter in this setting.}

\begin{maintheorem}
\label{t:gap_bounds}
With the same conditions as in Theorem \ref{T:HKLRP}, let $\mathcal{C}^1(N)$ denote the largest component of the graph induced by the vertices $[-N,N]^d$.  Then there exists a constant $\delta$ such that asymptotically almost surely
\[
\hbox{Gap}(\mathcal{C}^1(N)) \geq N^{d-s} \log^\delta N
\]
where $\hbox{Gap}$ denotes the spectral gap of the SRW on $\mathcal{C}^1(N)$.
\end{maintheorem}

Simple conductance arguments (e.g. \cite{BBY}) show that this is tight up to logarithmic factors.  This extends the results of \cite{BBY} in two directions: firstly to higher dimensions and secondly and more significantly to the setting of LRP where nearest neighbour bonds are not included with probability 1.

A main motivation for this work is to provide an important estimate for our companion paper \cite{CS} where we determine the scaling limits for random walks on supercritical LRP clusters $s \in (d, d+1)$ for $d \geq 1$.  The scaling limit is a L\'{e}vy process with exponent $\alpha=s-d$ which confirms a conjecture of \cite{Berger-Biskup}.  As mentioned above, in $d=1$ we also prove convergence to Brownian motion for $s > 2$ as well.  

%% maybe we shouldn't discuss the limitations of our later work until that paper %%% For $d \geq 2$, we do not address the case $s \in [d+1, d+2)$ as these involve more subtle issues as the limiting process is not a finite variation process.

%Our proof of this result naturally breaks into two distinct pieces.  An annealed limit theorem for SRW may be proved using methods from ergodic theory.  To extend the annealed limit theorem to a quenched limit law, we require some amount of control on the return probabilities $P_t^\omega(0, 0)$ for SRW started from the infinite component of the configuration.  Thus, Theorem \ref{T:HKLRP} provides the key heat kernel bounds required for SRW on LRP clusters.  

%The papers \cite{Kip-Var, Demasi, Sid-Szn, RM, Barlow-HKP, Mathieu-Piat, Berger-Biskup} consider aspects of random walks on super critical percolation clusters.  

%%%%%%  We maybe don't need that detailed an explanation in this paper%%%%%%We remark that our limit law proof is quite different from \cite{Sid-Szn, Mathieu-Piat, Berger-Biskup}; indeed, our setting represents more of a law of large numbers than a CLT.  Interestingly, those authors also required some \textit{ a priori} control of the heat kernel of SRW on the corresponding percolation cluster as well.  Fortunately, Barlow \cite{Barlow-HKP} had in the meantime derived tight bounds on the matrix elements of the heat kernel for SRW on a nearest neighbor supercritical percolation cluster, see also \cite{RM}.

It remains open to remove the logarithmic factors in the upper bound of Theorem \ref{T:HKLRP} which would be necessary to obtain sharp asymptotics analogous to \cite{Barlow-HKP}.  Our approach, following \cite{BBY}, of using the multi-commodity flows to bound the spectral gap of the SRW seems to inherently lose logarithmic factors.
%{Up to $\log$ corrections, we plan to address the question of matching heat kernel lower bounds in a future paper \cite{BCS}.  There we will also address the question of intersections of independent SRWs and the question of \cite{BBY} on the spanning forest on LRP.}  %%% let's see what the status is at the time of posting

\subsection{Proof Outline}
\label{S:impression}
The goal is to adapt, in so far as it is possible, the methods for obtaining heat kernel upper bounds employed in \cite{Barlow-HKP}, particularly Section 3.1.  This causes several complications: since the underlying percolation process is, in a sense, scale invariant, it is not at all clear how to choose balls on which we may obtain good spectral gap bounds.  Moreover, since the volume growth of the largest cluster is exponential in graph distance (see \cite{BiskupLR, 2009PT}) none of the standard methods apply to our setting.  

Our solution to these problems are twofold.  
We study SRW on the graph $G$, defined as the largest percolation cluster inside a ball of $\ell^{\infty}$ radius $N$ about the origin $B_N(0)$ where $N$ is chosen sufficiently large.  Note that because of the scale invariance of LRP, this choice of $N$ will be typically much larger than the time scale on which we wish to study the walk. 
While Barlow's argument was made for the particular case of nearest neighbour percolation, abstractly it only requires a partitioning of the percolation cluster into connected subsets, all of which are "good" in that they possess similar volumes, small diameters, and good spectral gap bounds.  His proof is inherently multi-scale and thus one must adjust the choice of partition on each time scale.  

For each time $s$ we construct a partition $\PP_s$ of the largest component of the large ball $B_N(0)$ into disjoint connected subsets which all have similar volumes, small diameters, and "good" spectral gap bounds.  As some vertices are only connected to the giant component through very long connections to distant parts of the graph this must be done carefully and can not involve only the local considerations.  We then adapt the arguments of Barlow replacing his ``good balls" by elements of the partition $\PP_s$ to obtain the stated heat kernel decay.

The issue becomes how to construct $\PP_s$ {(and this represents the bulk of the work)}.  The basic idea, which in itself is not sufficient, is to consider a re-normalized lattice, tiling $B_N(0)$ by blocks of side-length $\log^{\gamma} N$ for some appropriately chosen $\gamma$.  These blocks form the lattice points of the re-normalized graph.  We proceed by revealing edges in a number of stages. At each stage, the remaining unrevealed edges are independent of one another and of the construction to that point.

For the first stage, denote any "microscopic" block by $\Lambda$.  We reveal only the edges between vertices within $\Lambda$ and denote the largest component by $\CC^*(\Lambda)$.  Biskup \cite{BiskupLR} proved that with exponentially high probability, the size of $\CC^*(\Lambda)$ is at least $\rho \log^{d \gamma} N$ and we say the first step succeeds this estimate holds for all blocks.

The second stage is conditional on the success of the first stage: we reveal all edges between vertices within these largest components $\CC^*(\Lambda)$. On the event that the first stage succeeds for two neighbouring blocks $\Lambda, \Lambda^*$, the probability that they are directly connected is bounded below by
\be
\BbbP(\CC^*(\Lambda) \leftrightarrow \CC^*(\Lambda^*)) \geq 1- e^{-\beta \log^{(2d-s) \gamma} N }.
\ee
Thus we may guarantee all blocks are directly connected to their neighbors with high probability.  We say the second stage succeeds if the largest components of all neighbouring blocks are connected.  This cluster of vertices forms what we refer to as the \textit{core} for the maximal connected cluster of $B_N(0)$.  

On the event that the first two stages are successful, the core defines a re-normalized lattice isomorphic to $B_{N/\log^\gamma N}(0)$.  Moreover revealing the remaining edges of the core of $B_N(0)$ naturally induces a percolation process on the re-normalized graph, which dominates LRP with \textit{all} {nearest neighbor} connections present.  Furthermore, the effective value of $\beta$ for the induced LRP process \textit{increases} from $\beta$ to $\beta \log^{(2d-s)\gamma} N$ which results in a very dense well connected graph.

If we were studying the SRW on the graph induced by the core, we could apply the methods of \cite{BBYr} by using the components $\CC^* (\Lambda)$ as building blocks.
However, at this point we are confronted with another problem: because we do not assume nearest neighbor connections occur with probability one, after the first two steps there are many small clusters left over which are nonetheless part of the largest component, possibly through long connections.
The main hurdle to overcome is how to allocate small clusters to parts of the core so that our diameter, volume and spectral gap bounds still hold.  This is done through a multistep clustering scheme, which we discuss in detail below.  Again we reveal edges in stages and take advantage of the independence in LRP.

The rest of our paper is organized roughly in reverse order of this discussion:  In the next section we study the construction of the core vertices and then demonstrate how to allocate small clusters to the core so as to maintain regularity properties.  Then we prove spectral gap bounds for SRW on these allocated clusters. 
In Section \ref{S:HK} we adapt these probabilistic estimates to obtain the desired heat kernel bounds.

\section{Graph structure}
The construction of the long range percolation measure is standard, but we recapitulate some aspects here to collect notation.  We begin by introducing $\Lambda_N= [-N, N]^d \cap \mathbb Z^d$.  Given a subset $\Lambda \subset \mathbb Z^d$, let $\EE_{\Lambda}$ denote the set of unordered pairs $\la x, y\rra$ such that one of $x, y \in \Lambda$.  If $\Lambda= \mathbb Z^d$, let $\EE= \EE_{\Z^d}$.  An edge $b \in \EE$ is said to be occupied (relative to $\omega \in \Omega$) if $\omega_b=1$.  By slight abuse, we shall regard $\omega$ as a random graph with vertices $\Z^d$ and edges given by $\{b: \omega_{b}=1\}$. 

%Let $\FF_{\Lambda}$ denote the usual Borel $ \sigma$-field on $\{0,1\}^{\EE_{\Lambda}}$ associated with the product topology and let $\mu_N = \prod_{\la i j\rra\in \EE_{\Lambda_N} }p_{ij}$, where the $p_{ij}$ denote translation invariant independent Bernoulli edge distributions of the form \eqref{Eq:Cond}.  Here $\|i-j\|_2$ denotes the $\ell_2$ distance.  We assume the $p_{ij}$ are \textit{percolating}. 

%{\blue Let $\Omega_{f,N}=\{\omega \in \{0,1\}^{\EE_{\Lambda_N}}: \#\{b :\omega_b =1\}< \infty\}$. }
%A: Since this is a finite graph won't the degrees always be finite?
%N: No you misread my definition of \EE_{\lambda};  it IS inifinite
%Then since $s \in (d, d+2)$, $\mu_N$ is in fact supported on $\Omega_{f, N}$.
%The family $(\mu_N, \Omega_{f, N}, \FF_{\Lambda_N})$ form a consistent family of probability distributions, so the Kolmogorov Extension Theorem implies the existence of a unique limiting measure $\mu$ on $\Omega = \{0,1\}^{\EE}$ with the standard product $\sigma$-field.  

Recall that, given a connected finite graph $G=(V, \EE)$, SRW on $G$ defines a transition kernel $P_{x, y}$, which as an operator is self adjoint with respect to $\ell^{2}(\pi)$ where 
\be
\pi(x)\propto \textrm{deg}_{G}(x)
\ee
is the stationary probability measure for $P$.  The spectral gap is defined as
\be
\textrm{Gap}_{G} := \inf_{f \perp \mvone} \frac{(f, (\textrm{Id}-P)f)_{\pi}}{(f, f)_{\pi}}
\ee

For a connected subset $\CC \subset \omega$, it will at times be useful to consider the subgraph induced by $\CC$ from $\omega$ and the SRW on this subgraph.  We will denote the transition kernels for these walks by
\be
P^{\omega}_{t}(x, y) \text{ and, respectively } P^{\omega}_{\CC, t}(x, y)
\ee
as necessary.

Throughout, we work with the $\ell^\infty$ metric on $\Z^d$ and various intrinisic graph metrics on percolation clusters $\CC$.  We will generally be careful to note which of these metrics we are using by either an explicit statement or by the notation
\be
d_{\Z^d, \infty}(\cdot, \cdot) \quad d_{\CC}(\cdot, \cdot)
\ee
respectively.  Note that for any $\Z^d$, $d_{\Z^d, \infty}(\cdot, \cdot) $ is comparable to $
d_{\Z^d, 2}(\cdot, \cdot) $, the usual Euclidean metric.
Let $B_{N}(0)$ be the $\ell^{\infty}$ ball about $0$ or radius $N$.  Given $N_1< N$, let $\scrO=\scrO(N, N_1)$ denote a fixed minimal cover of $B_{N}(0)$ by boxes $\Lambda_{N_1}$ of the form $\Lambda_{N_1}= x+ B_{N_1/2}(0)$ so that $|\Lambda_{N_1}| = N_1^d$.  Note that we implicitly assumed $N_1$ is even and also note also that if $N_1$ does not divide $N$, this cover may not be unique; neither will be of any loss and we fix some deterministic rule for choosing a cover.  Let $\PP_{N_1}$ denote the collection of vertex disjoint boxes $\Lambda_{N_1}$ whose union is the cover $\scrO$.

\begin{theorem}
\label{T:gap}
Suppose that $p_{x, y}$ satisfy the conditions surrounding \eqref{Eq:Cond}.  Suppose that $d<s<(d+2) \wedge 2d$.  Then there exist universal constants $\rho_1, \rho_2, C_1, C_2, c, C, \delta_1, \delta_2>0$ so that if $N\geq N_1$ and if $\scrO$ is a minimal cover of $B_N(0)$ by boxes of side-length $N_1$ then there is an event $A(N, N_1)$ with
\be
\mu(A(N, N_1)) \geq1- C_1 e^{-C_{2} \log^2 N_1},
\ee
and on $A(N, N_1)$, we can find a partition $\PP$ of $\scrO$ 
%A:Should this be B_N(0)?
%N:No, I mean to partition \scrO
into connected subsets $\left(\scrG(0, \Lambda_{N_1})\right)_{\Lambda_{N_1} \in \PP_{N_1}}$ so that:
\begin{enumerate}
\item $Vol(\left(\scrG(0, \Lambda_{N_1})\right) \cap \Lambda_{N_1}) \geq \rho_1 N_1^d$.

\item $Vol\left(\scrG(0, \Lambda_{N_1})\right) \leq \rho_2 N_1^d$.

\item $Diam (\scrG(0, \Lambda_{N_1}) \leq c \log^{\delta_1}(N_1)$.

\item $\textrm{Gap}_{\scrG(0, \Lambda_{N_1})} > C N^{d-s}_1\log^{-\delta_2} N \textrm{ for all } \Lambda_{N_1} \in \PP_{N_1}$.
\end{enumerate}
where the diameter is in the graph distance on the induced components.
\end{theorem}

Let us define $\tau_x \cdot \omega$ to be the usually shift of the conifiguration $\omega$:
\begin{equation}
\tau_x \cdot \omega_{\la y, z \rra}=  \omega_{\la x+y, x+z \rra}
\end{equation}

Let $\tau_x \cdot A(N, N_1)$ denote the image of $A(N, N_1)$ under the shift $\tau_x$:
\begin{equation}
\tau_x \cdot A(N, N_1)= \{ \omega: \tau_x\cdot \omega \in A(N, N_1)\}.
\end{equation}
Let us define $B(\epsilon, N, x) = \cap_{N^{1-\epsilon} \geq N_1 \geq N^{\epsilon} } \tau_x \cdot A(N, N_1)$.  %Not quite sure the role $x$ plays here?
Then it is easy to see:
\begin{corollary}
\label{C:Gap}
For all $x\in \Z^d$ and all $\epsilon>0$, there exists an integer valued random variable $N_{\epsilon, x}(\omega)$, $N_{\epsilon, x}(\omega)< \infty$ $\mu\: a.s.$ such that  for all $N>N_{\epsilon, x}(\omega)$,
$B(\epsilon, N, x)$ occurs.
Moreover the tail of $N_{\epsilon, x}(\omega)$ satisfies
\be
\mu(N_{\epsilon, x}(\omega)>k)< e^{-C_3 \log^2 k}
\ee
for some universal constant $C_3$.
\end{corollary}

The proof of this theorem is a bit involved, taking the better part of the next $15$ pages.  We already outlined the issues encountered in Section \ref{S:impression}.
We begin by building, with high probability, a \textit{core} of vertices, which may be partitioned into connected subsets of $\scrO(N, N_1)$ which have good diameter and volume growth in their own right.  The second step is to argue that we may allocate vertices \textit{not} in the core to connected subsets of the core in such a way that the diameters and volumes do not change in appreciably. By construction we will be able to obtain spectral gap on the partition sets by using a multi-commodity flow argument \cite{Sinclair} and comparison to the corresponding result for the Erd\"{o}s-R\'{e}nyi graph with connection probability in the connected regime.

\subsection{Step 1: The Core}
Let us denote by $\CC^*(\Lambda_L)$ the largest component of $\Lambda_L$ \textit{connected inside $\Lambda_L$}.  
Let us recall a crucial result from \cite{BiskupLR}:
\begin{theorem}
\label{thm3.2}
Let~$d\ge1$ and consider the probabilities~$(p_{xy})_{x,y\in\Z^d}$ such that \eqref{Eq:Cond} holds for some~$s\in(d,2d)$. Suppose that~$(p_{xy})_{x,y\in\Z^d}$ are percolating. For each $s'\in(s,2d)$ there exist numbers~$\rho>0$ and~$L_0<\infty$ such that for each~$L\ge L_0$,
\begin{equation}
\label{2.8}
\mu(|\CC^*(\Lambda_L)|<\rho|\Lambda_L|\bigr)\le e^{-\rho L^{2d-s'}}.
\end{equation}
\end{theorem}

%By Theorem \ref{thm3.2}
%\be
%\label{Eq:LL}
%\mu(\exists \Lambda_L \in\PP_L : |\scrC^*(\Lambda_L)|<\rho|\Lambda_L|) \leq \left(\frac{N}{L}\right)^d e^{-\rho L^{2d-s'}}
%\ee
In particular, once~$L$ is sufficiently large, the largest connected component in~$\Lambda_L$ typically contains a positive fraction of all sites in~$\Lambda_L$.
%Moved out of the theorem
We will only invoke this theorem once.  Fix, once and for all, $s', \rho$ and $L_0$ for which the theorem holds.
We choose a set of scales on which our construction will take place:
With $N, N_1$ fixed, we choose (if possible) 
\be
\label{eq:chooseN} 
N= N_0 > N_{1}> N_2 > N_3>N_4
\ee 
so that 
\be
N_{2} = \lfloor N_{1}^{(s-d)/d}\log ^{3/d} N_1 \rfloor, N_3=\lfloor N_2^{1/2}\rfloor\text{ and }N_4= \lceil\log^{2/(2d-s')} N\rceil.
\ee  
We remark that, for our objectives, the choice of $N_2$ is necessary, $N_3$ has a lot of choice, and $N_4$ is roughly necessary.  For convenience we will suppose $N_{j}$ divides $N_{j-1}$ for $j \geq 2$, though this is patently false.  It is quite standard to modify the proofs to take into account this discrepency (for example by working with a dyadic decomposition).

For any box $\Lambda_L$ and for any integer $K$ dividing $L$, let $\tilde P_K( \Lambda_L)$ denote the partition of $\Lambda_L$ into boxes of side length $K$.
Let $(\PP_{N_j})_{1 \leq j \leq 4}$ denote a sequence of refinements of $\scrO$ defined as follows.  Let $\PP_{N_1}= \{\Lambda_{N_1}: \Lambda_{N_1} \subset \scrO\}$.  Inductively, having chosen the refinement $\PP_{N_j}$ to consist of boxes with side length $N_j$, let $\Lambda_{N_j} \in \PP_{N_j}$.  As above we may introduce $\tilde \PP_{N_{j+1}}(\Lambda_{N_j})$, a partition of $\Lambda_{N_j}$ into blocks of side-length $N_{j+1}$ and then let
\be
\PP_{N_{j+1}}:= \cup_{\Lambda_{N_j} \in \PP_{N_j}} \tilde \PP_{N_{j+1}}(\Lambda_{N_j})
\ee

\begin{definition} 
We shall say that two blocks $\Lambda_1, \Lambda_2 \in  \PP_{N_{j+1}}$ are adjacent
if $d_{\Z^d, \infty}(\Lambda_1, \Lambda_2) \leq 2$, where $d_{\Z^d, \infty}$ is the Hausdorff distance measured with respect to the $\ell^{\infty}$ norm on $\Z^d$.
\end{definition}

For each $\Lambda_{N_j}$, we construct a random set
$\scrG(\Lambda_{N_j})\subset \Lambda_{N_j}$, as follows.
We define the cores intersection with a block $\Lambda_{N_4}\in \PP_{N_4}$ as 
\be
\scrG(\Lambda_{N_4})=\CC^*(\Lambda_{N_4}).
\ee
and say that $\Lambda_{N_4}$ is occupied if $| \scrG(\Lambda_{N_4})| >\rho|\Lambda_{N_4}|$.
%I think it might be cleaner to define \scrG(\Lambda_{N_4}) always but say it is occupied if it is big enough

We use interchangeably the terminology that $\scrG(\Lambda_{N_4})$ exists, has been constructed or that $\Lambda_{N_4}$ is occupied.  Let $\FF_{*}$ denote the $\sigma$-algebra generated by all edge events $\omega_{x,y}$ for pairs of vertices $x,y\in \scrO$ in the same box in the partition $\PP_{N_4}$ and let $\OO\in \FF_{*}$ denote the event that all $\scrG(\Lambda_{N_{4}})$ are occupied.  
By Theorem \ref{thm3.2} and a union bound we have that
\begin{equation}\label{e:prob_OO}
\mu(\OO) \geq 1 - e^{-c\log^2 N}
\end{equation}

Let $\FF_j$ denote the $\sigma$-algebra
\[
\FF_j =\sigma\left\{\FF_* \vee \sigma\{ \omega_{x,y}:  |x-y|<N_j\in \Lambda_{N_j}\in \PP_j, x,y \in \cup_{\Lambda_{N_4} \subset \Lambda_{N_j}} \scrG(\Lambda_{N_4}) \} \right\}
\]
Let $\{\scrG(\Lambda_{N_{4}}) \leftrightarrow \scrG(\Lambda_{N_{4}}')\}$ denote the event that there is a direct connection from $\scrG(\Lambda_{N_{4}})$ to $\scrG(\Lambda_{N_{4}}')$ and $\{\scrG(\Lambda_{N_{4}}) \nleftrightarrow \scrG(\Lambda_{N_{4}}')\}$ denotes the complimentary event.  Let $\AA$ denote the event that $\{\scrG(\Lambda_{N_{4}}) \leftrightarrow \scrG(\Lambda_{N_{4}}')\}$ occurs for each pair $\Lambda_{N_4}, \Lambda_{N_4}'$ of adjacent blocks.

%{\magenta Next, given $\Lambda_{N_j}$, assume that random sets
%\be
%\{\scrG(\Lambda_{N_{4}}): \Lambda_{N_{4}} \subset \Lambda_{N_j}\}
%\ee
% have been constructed successfully.   We may next sample all edges in $\Lambda_{N_{j}}$ \textit{only} between the random sets 
% \be
% \{\scrG(\Lambda_{N_{4}}): \Lambda_{N_{4}} \subset \Lambda_{N_j}\}
% \ee 
% in \textit{adjacent} blocks $ \Lambda_{N_{4}},  \Lambda_{N_{4}}'$. }\

\begin{definition}
Let us define the \textbf{core} to be the graph $\scrG=\scrG_N$ with vertex set $\cup_{\Lambda_{N_4} \in \PP_{N_4}} \scrG(\Lambda_{N_4})$ and edges given by revealing all edges of $\ell^\infty$ length at most $N_2$ between these vertices. 
Similarly, sampling edges \textit{ inside } $\cup_{\Lambda_{N_4} \subset \Lambda_{N_j}} \scrG(\Lambda_{N_4})$ of length at most $N_2 \wedge N_j$, call the derived %it's not quite an induced graph since we don't take all the edges
 graph $\scrG(\Lambda_{N_j})$. 
\end{definition}

%Let us record useful properties regarding the $\scrG$ {\green and} $\scrG(\Lambda_{N_j})$.  Let $\FF_{N_j}$ denote the filtration defined by sampling all edges in $\scrC$ of (Euclidean) length less than $N_j$.
%Let
%\be
%\sigma(\{\scrG(N_{4})\}: \Lambda_{4} \subset \Lambda_j\})
%\ee
%denote the $\sigma$-algebra generated by the construction of the family
%\be
%\{\scrG(N_{4})\}: \Lambda_{4} \subset \Lambda_j\}
%\ee
%Conditionally on the existence of $\scrG(\Lambda_{N_{4}}), \scrG(\Lambda'_{N_{4}})$, let $\sigma( \EE_{\scrG(\Lambda_{N_{4}}), \scrG(\Lambda'_{N_{4}})})$ denote the $\sigma$-algebra generated by the variables
%\be
%\{\omega_{x, y}: x \in \scrG(\Lambda_{N_{4}}), y \in \scrG(\Lambda_{N_{4}}'), \|x-y\|_{\infty}< N_2\}.
%\ee
%Let
%\be
%\sigma(\scrG(N_j))= \sigma\left(\FF_{N_4}, \sigma( \EE_{\scrG(\Lambda_{N_{4}}), \scrG(\Lambda'_{N_{4}})}: \Lambda_{N_{4}}, \Lambda'_{N_{4}} \subset \Lambda_{N_j})\right).
%\ee
%Properly, we should enlarge our $\sigma$-algebras to include the situation that $\scrG(\Lambda_{N_4})$ does not exist, but we will always work on the event that all do so this is superfluous in this case. 
\begin{lemma}
\label{L:CoreProps}
There exist constants $\delta_1, c_1, c_2 > 0$ such that for all $N> N_1>N_2> N_3>N_4$ chosen as in \eqref{eq:chooseN},  
there exists an event 
\be
\textrm{Co}(N, N_1)= \textrm{Co}(N, N_1, N_2, N_3, N_4)
\ee
with
\be
\mu(\textrm{Co}(N, N_1)) \geq 1- c_1e^{- c_2\log^2 N}
\ee
satisfying the following properties:
\begin{itemize}
\item\textbf{Connectedness:} $\scrG(\Lambda_{N_j})$ is connected \textit{in} $\Lambda_{N_j}$.  If  $d_{\Z^d, \infty}(\scrG(\Lambda_{N_{4}}), \scrG(\Lambda_{N_{4}}' ))< N_4$ then $\scrG(\Lambda_{N_{4}}), \scrG(\Lambda_{N_{4}}')$ are directly connected. 
\item \textbf{Volume Growth:} $|\scrG(\Lambda_{N_j})| \geq \rho |\Lambda_{N_j}|$.
\item  \textbf{Diameter Bounds:}  $\textrm{ Diam }(\scrG(\Lambda_{N_j})) <  N_4^d \log^{\delta_1} N_j$ for all $\Lambda_{N_j}$ for $j = 2, 3$.
\end{itemize}
Also, we have
\begin{itemize}
\item{\textbf{Localization:}} $\scrG(\Lambda_{N_j}) \subset \Lambda_{N_j}$.
\item
\textbf{Measurability:} $\scrG(\Lambda_{N_j})$ is measurable with respect to the $\sigma$-algebra $\FF_{2\vee j}$.  In particular the construction of $\scrG$ does not call on edges longer than $N_2$.
\end{itemize}

\end{lemma}
We remark that the diameter bounds here are not optimal.  Bounds of the correct order have been obtained in \cite{BiskupLR}.  By our choice of the $N_j$, we could have absorbed $N_4^d$ into $\log^{\delta_1} N_j$, but we keep in this form to track the dependence of estimates on our choices.

%Given $\Lambda_{N_j}$, define the event that 
%\be
%\OO_{\Lambda_{N_j}} =\{\scrG(\Lambda_{N_4}): \Lambda_{N_4} \subset \Lambda_{N_j} \textrm{ are occupied and adjacently connected }\}.
%\ee 

The key observation is that if we regard $
\{\scrG(\Lambda_{N_{4}}):  \Lambda_{N_{4}} \subset \Lambda_{N_j}\}$
as vertices of a course grained lattice in a box $B_{N_j/N_4}$, then on the event $\OO \cap \AA$, sampling all edges between all remaining vertices of $
\{\scrG(\Lambda_{N_{4}}): \Lambda_{N_4} \subset \Lambda_{N_j}\}$ stochastically 
dominates long range percolation on a box of side length $N_j/N_4$ \textit{ with all nearest neighbor edges present}.  
Indeed, if 
\be
d_{\Z^d, \infty}(\scrG(\Lambda_{N_{4}}), \scrG(\Lambda_{N_{4}}')) \leq \ell N_4
\ee 
then
\be
\label{Eq:Key2}
\mu(\scrG(\Lambda_{N_{4}}) \nleftrightarrow \scrG(\Lambda_{N_{4}}') \mid \OO) \leq e^{-2^{d/2}\rho^2 (\ell+1)^{-s} N_4^{2d-s}}
\ee
First we note that with a union bound this bound implies that
\begin{equation}\label{e:prob_AA}
\mu(\AA|\OO) \geq 1 - e^{-c \log^2 N}
\end{equation}
Further on $\OO \cap \AA$, the re-normalized block percolation process determined by the edges between these $\scrG(\Lambda_{N_4})$ stochastically dominates long range percolation with exponent $s$ and $\beta=\beta_{N_4}=N_4^{2d-s}$ \textit{where nearest neighbor blocks are conditioned to be connected with probability one}.  Note that by our definition of $N_4$, this leads to a highly connected graph.

Note that $\scrG(\Lambda_{N_j}) \subset \scrG$ for all $\Lambda_{N_j} \in \PP_{N_j}, j=1, 2, 3, 4$.
Before giving the proof (most of which is clear by construction) we record a diameter bound separately for easy reference.

%Let 
%\be
%\OO:= \{\text{ $\Lambda_{N_4}:\Lambda_{N_4} \in \PP_{N_4}, $ are occupied }\},
%\ee
%and
%\be
%\AA:= \{ \scrS(\Lambda_{N_4}) \text{ directly connected to } \scrS(\Lambda_{N_4}'):d_{\Z^d, \infty}(\Lambda_{N_4}, \Lambda_{N_4}') \leq 1 \}
%\ee

\begin{lemma}
\label{L:Diam}
Consider the family $(\scrG(N_j)_{\Lambda_{N_j}}$ on the event $\OO \cap \AA$.
Then
\be
\mu(\exists \Lambda_{N_j}: \textrm{ Diam }(\scrG(\Lambda_{N_j})) > \log^{\delta} N_j N_4^d \textrm{ some } \Lambda_{N_j}| \OO \cap \AA) \leq 2^dN^d e^{-N_4^{2d-s}}
\ee	
for $j =2, 3$.  %Since we don't include long edges
\end{lemma}
\begin{proof}
By \eqref{Eq:Key2}, we may apply:
\begin{corollary}[Corollary 5.1 from \cite{CGS}]
Let $B_K$ be a box of side length $K$.  and consider long range percolation on $B_K$ with exponent $s$ and parameter $\beta$ and nearest neighbors connected $a.s.$  For any $C>0$ there exists $\delta$ such that
\be
\mu(\textrm{Diam }(B_K)> \log^\delta K) \leq e^{-\Theta(\log^C K)}.
\ee
\end{corollary}
%  I don't really understand how we use this Corollary because can't we take C=2 and get the result straight away?
We amplify this estimate via the observation from \eqref{Eq:Key2} that the course grained process with vertices $(\scrG(\Lambda_{N_4}))_{\Lambda_{N_4} \subset \Lambda_{N_j}}$ dominates a long range percolation process, denoted by $\LL \RR$, with exponent $s$ and $\beta_{N_4}= \beta \rho^2 N_4^{2d-s}$.

The point is that we can instead view $\LL \RR$ as a union of $\lfloor \rho^2 N_4^{2d-s} \rfloor$ independent identically distributed long range processes $(\LL \RR_{j})_{j=1}^{\lfloor \rho^2 N_4^{2d-s} \rfloor}$ with exponent $s$ and parameter $\beta' \geq \beta$.  By slight abuse of notation, we let $\mu$ denote the coupling measure between the process on $(\scrG(\Lambda_{N_4}))_{\Lambda_{N_4} \subset \Lambda_{N_j}}$ and the family  $(\LL \RR_{j})_{j=1}^{\lfloor \rho^2 N_4^{2d-s} \rfloor}$.

Thus applying the Corollary to each independent copy, we have
\be
\mu(\nexists j : \textrm{ Diam }(\LL \RR_{j}) < \log^\delta N_j N_4^d \textrm{ some } \Lambda_{N_j}| \OO \cap \AA) \leq e^{-\lfloor \rho^2 N_4^{2d-s}\rfloor \Theta(\log^C N_j)}
\ee
Taking a union bound finishes the proof.
\end{proof}

\begin{proof}[Proof of Lemma \ref{L:CoreProps}]
The volume growth, connectedness, localization and measurability claims all follow by construction, on the event that $\OO\cap \AA$.  By equations \eqref{e:prob_OO} and \eqref{e:prob_AA} this event holds with probability at least $1- e^{- c \log^2 N}$.  The diameter bound holds by applying Lemma \ref{L:Diam} and taking a union bound.

%{\magenta For the remaining claims, let us begin by computing the probability that \textit{all} $\Lambda_{N_4}$ are occupied and all pairs $\scrG(\Lambda_{N_4}), \scrG(\Lambda_{N_4}')$ are directly connected if $\Lambda_{N_4}, \Lambda_{N_4}'$ separated by at most $2$ in $\ell^\infty$.

%By Theorem \ref{thm3.2},
%\be
%\mu(\exists \Lambda_{N_4} \in \PP_{N_4}, \text{ $\Lambda_{N_4}$ not occupied }) \leq 2^dN^d e^{-\rho^2 N_4^{2d-s'}}
%\ee

%On the event 
%\be
%\OO:= \{\text{ $\Lambda_{N_4}:\Lambda_{N_4} \in \PP_{N_4}, $ are occupied }\},
%\ee
%the event 
%\be
%\AA:= \{ \scrG(\Lambda_{N_4}) \text{ directly connected to } \scrG(\Lambda_{N_4}'):d_{\Z^d, \infty}(\Lambda_{N_4}, \Lambda_{N_4}') \leq 2 \}
%\ee
%satisfies
%\be
%\mu(\AA^c|\OO) \leq 6^d N^d e^{-2^{-s}\rho^2N_4^{2d-s}}.
%\ee

%Furthermore, on $\OO$, $\textrm{Diam }{\scrG(\Lambda_{N_4}}) \leq N_4^d$.  By Lemma
%\ref{L:Diam}, for any $j \leq 3$,
%\be
%\mu(\exists \Lambda_{N_j}: \textrm{Diam }(\scrG(\Lambda_{N_j})> N_4^d \log^{\delta_1} N_j|\OO\cap\AA) \leq 2^dN^d e^{-N_4^{2d-s}}.
%\ee
%Taking intersections, collecting estimates and recalling the definition of $N_4$ finishes the proof.}\
\end{proof}

\subsection{Step 2:  Allocating Small Clusters to the Core}

Recall that we are aiming to show that we can, with high probability, partition the largest cluster of a minimal cover $\scrO$ of $B_N(0)$ into small subsets, each of which has similar volume and spectral gap.  This has already been partially achieved for the core $\scrG$.  Significantly for our future bounds on spectral gaps, $\scrG$ was constructed without knowledge of edges longer than $N_2$.

This core is quite large, and the next lemma shows (among other things) that with high probability it is in the largest cluster.  Our partition, however, must include \textit{all} vertices in the largest cluster, not just those in the core.  We proceed to allocate the remaining vertices not in the core to the subsets $\scrG(\Lambda_{N_1})$ so as to not alter the volume, diameter and (therefore) spectral gaps of the $\scrG( \Lambda_{N_1})$ already constructed by too much (Note: we have not yet addressed the spectral gaps $\scrG( \Lambda_{N_1})$, we will directly analyze the gap of our full partition).
The following lemma justifies our analysis of the core, showing that with high probability it is part of the largest component, and furher serves as a model for future development:
\begin{lemma}
\label{L:Sizes}
Let $n_1 > n_2 \geq \dotsc \geq n_m $ enumerate the cluster sizes inside $\scrO$ (having sampled all internal edges).  Then there exists $c_3, c_4,\delta_2> 0$ such that
Then
\be
\mu(n_2 > N_4^{s-d} \log^{2} N) < c_3e^{-c_4\log^2 N}.
\ee
Moreover, the core $\scrG$ is a subset of the largest component except with probability $O(e^{-c_4\log^2 N})$. 
\end{lemma}

\begin{proof}
On the event $\OO$, we next reveal the edges connecting all vertices in $\scrG^c\cap \scrO$ and denote the $\sigma$-algebra generated by these edges as $\BB$.  The set $\scrG^c\cap \scrO$ is broken into a sequence of clusters 
\be
\CC^1, \CC^2, \dotsc, \CC^M \subset \scrG^c\cap \scrO
\ee 
where the edges between $\scrG$ and $(\CC^1, \dotsc \CC^M)$ have not yet been revealed unless they were revealed in $\FF_*$, that is pairs of vertices which are in the same box in $\PP_{N_4}$.  On $\OO$ the core $\scrG$ has positive density down to the level $N_4$.  Thus
\be
\mu(\CC^i \nleftrightarrow \scrG |\BB,  \FF_*, \OO ) \leq e^{-  2 \rho \beta N_4^{d-s} |\CC^i|}
\ee
since each for vertex $x\in \CC^i$ there are at least $\rho N^d_4$ vertices in the core in an adjacent block in $\PP_4$ and $x$ could be connected to them each independently with probability at least $1-e^{-  2 \beta N_4^{-s} |\CC^i|}$.  As the total number of clusters is bounded by $2^dN^d$, by a union bound except with probability $2^d N^d e^{-  2 \beta N_4^{-s} |\CC^i|}$ there is no component of size greater than $N_4^{s-d} \log^{2} N$ which is not connected to the core.

On the event $\OO\cap \AA$, which occurs with probability at least $1-e^{-c\log^2n}$,  the core forms a single component of size at least $\rho N^d$.  Hence the core is in the largest component and the second largest component is of size at most $N_4^{s-d} \log^{2} N$ except with probability  $c_3e^{-c_4\log^2 N}$ as required.

\end{proof}

To get around issues outlined above, we will reveal edges with endpoints in $\scrG^c\cap \scrO$ in stages.  Let us introduce notation to facilitate this step.  % {\magenta Let $\sigma(\scrG)$ denote the $\sigma$-algebra defined by the construction of $\scrG$,} %%%%%%%I don't think we use this notation???
Let $\HH_{N_j}$ denote the filtration defined by 
\[
\HH_{N_j} =\sigma\left\{\FF_2 \vee \sigma\{ \omega_{x,y}:  |x-y|_\infty <N_j, x,y \in \scrO \} \right\}.
\]
that is $\FF_2$ along with  all edges in $\scrO$ of $\ell^{\infty}$ length less than $N_j$ (only a nontrivial construction if $j=3,4$.  Let $\scrX_{N_j}$ denote the graph induced by $\scrG$ and sampling \textit{all} edges in $\scrO$ of length at most $N_{j}$  (thus $\scrX_{N_{j}}$ is measurable with respect to $\HH_{N_j}$).   

Let $\CC(j, \scrG)$ denote the connected component of $\scrG$ in this graph.  For $x \in \CC(j, \scrG)\backslash \scrG$ we define an allocation $\phi_{j}(x)$ (formally a function from $\CC(j, \scrG)$ to $\PP_{N_{2}}$).  Intuitively, we would like
\be
\phi_{j}(x)= \Lambda_{N_{2}} \text{ if }d_{\scrX_{N_j}}(x, \scrG \cap \Lambda_{N_2})=d_{\scrX_{N_j}}(x, \scrG)
\ee 
however this definition has two problems:
First, there may be ties and so strictly speaking the intuitive allocation is not well defined;
Second, it is important that $\phi_{j}^{-1}(\Lambda_{N_{2}})$ is a connected subset of the maximal component. 

Thus, to define $\phi_{j}$ on $\CC(j, \scrG)$ we proceed inductively. Consider the sets 
\be
R_k(\GG)=\{x\in \CC(j, \scrG) : d_{\scrX_{N_j}}(x,\GG)\leq k\}.
\ee
Here $d_{\scrX_{N_j}}(\cdot, \cdot)$ denotes the \textit{$\scrX_{N_j}$-graph} distance.
Now let $A_1(\Lambda_{N_2})=\{x\in R_1(\GG): d_{\scrX_{N_j}}(x,\GG\cap \Lambda_{N_2})\leq 1\}$
with ties broken according to some deterministic rule.  Inductively, suppose that $A_k(\Lambda_{N_2})$ have been defined and let
\be
A_{k+1}(\Lambda_{N_2})=\{x\in R_{k+1}(\GG) \backslash R_k(\GG) : d_{\scrX_{N_j}}(x,A_k(\Lambda_{N_2}))\leq 1\}
\ee
again with ties broken according to some deterministic rule.
Finally, let
\be
\phi_{j}(x):= \Lambda_{N_{2}} \text{ if }x \in \cup_k A_k(\Lambda_{N_2}).
\ee 
Then if
\be
\label{Eq:Part}
\scrG(j, \Lambda_{N_2})=(\scrG\cap\Lambda_{N_{2}}) \cup \{x: \phi_{j}(x)= \Lambda_{N_2}\},
\ee  
$(\scrG(j, \Lambda_{N_2}))_{\Lambda_{N_2} \subset \scrO}$ defines a partition of $\CC(j, \scrG)$ into connected subsets.

\subsubsection[Allocation Phase 1]

We shall need a simple bound on the maximal degree of a vertex in $B_K(0)$. Let $\textrm{deg}^{\omega}(x)= \sum_{y \in \Z^d} \mvone_{\omega_{x, y}=1}$ and $\textrm{Deg}(\Lambda_K)= \max_{x \in \Lambda_K} \textrm{deg}^{\omega}(x)$.  
\begin{lemma}
\label{L:Deg}
Then $\forall \delta > 0, \exists C(\delta)>0$ such that for all $\ell>0$, we have
\be
\mu(\textrm{Deg}(\Lambda_K)> b \log K) \leq C(\delta)K^d e^{- \delta b \log K}.
\ee
\end{lemma}
\begin{proof}
Using the exponential Markov inequality and independence of the summands of $\textrm{deg}^{\omega}(0)$,
\be
\mu(d^{\omega}(0) > \log K) \leq e^{-\delta\log K} \E_{\mu}[e^{\delta d^{\omega}(0)}].
\ee

But 
\be
\log \E_{\mu}[e^{\delta d^{\omega}(0)}] \leq C_{L} + \sum_{\|x\| \geq L} (e^{\delta}-1)\|x\|^{-s} \leq C'_L(\delta) + (e^{\delta}-1) L^{d-s}.
\ee
Taking a union bound finishes the proof.
\end{proof}

\begin{lemma}
\label{L:P1}
There exists $\rho_1, \rho_2, c_5, c_6, \delta_3>0$ such that for any $N> N_1> N_2>N_3> N_4$ as in \eqref{eq:chooseN} the following holds:
%for any $\delta\geq \delta_3$ and any constants $M_1, M_2, C_2, c_3, C_3$ the following holds.
Consider the graph $\scrX_{N_3}$.
We can find an event $E \subset \textrm{Co}(N, N_1)$ with $E \in \HH_{N_3}$ such that if $(\scrG(3, \Lambda_{N_2}))$ denotes the partition of $\CC(3, \scrG)$ into connected subsets from \eqref{Eq:Part} then on $E$
\be
\textrm{Diam}(\scrG(3, \Lambda_{N_2})) \leq C_2 \log^{\delta} N_2+ N_4^{s-d}\log^{2}N
\ee
and 
\be
 \rho_1 N_2^d \leq \textrm{Vol}(\scrG(3, \Lambda_{N_2})) \leq \rho_2 N_2^d
\ee
for all $\Lambda_{N_2}$.
Also, for any component $\CC \subset \scrX_{N_3}$ so that $\CC \neq \CC(3, \scrG)$, 
\be
Vol(\CC) \leq N_4^{s-d} \log^2 N.
\ee
Finally
\be
\mu(E|  \FF_{4} )\mvone_{\textrm{Co}(N, N_1)} \geq 1- c_5 e^{-c_6 \log^2 N} \mvone_{\textrm{Co}(N, N_1)}
\ee
$\mu$ a.s.
\end{lemma}
\begin{proof}
Let $\mu_{N_4}$ denote the conditional measure induced by $\mu(\cdot | \textrm{Co}(N, N_1), \HH_{N_4})$.  
Consider the law of $\scrX_{N_3}$ under $\mu_{N_4}$ and recall the definition of $\scrG(3,\Lambda_{N_2})$ \eqref{Eq:Part}. Let $B_{\ell}(x)$ denote the ball of $x$ with radius $\ell$ \textit{ measured in the graph distance on $\scrX_{N_3}$}.
Let $L(x)$ denote the least $\ell$ such that $|B_\ell(x)| \geq N_{4}^{s-d} \log^2 N$.  Then the same argument as in Lemma \ref{L:Sizes}
\be
\label{eq:a}
\mu_{N_{4}}(B_{L(x)}(x) \cap \scrG = \varnothing, L(x)< \infty) \leq e^{-c\log^2 N}
\ee
$\mu$ a.s.
and so
\be
\label{eq:b}
\mu_{N_{4}}(\exists x \text{ such that } B_{L(x)}(x) \cap \scrG = \varnothing, L(x)< \infty) \leq 2^d N^d e^{- c \log^2 N}
\ee
$\mu$ a.s.
Thus with high probability, $N_{3}$-clusters which exceed $N_{4}^{s-d} \log^2 N$ in volume  are directly connected to the core.
Since $L(x)$ is bounded by $N_{4}^{s-d} \log^2 N$ on the event that $L(x)< \infty$, we have good diameter bounds for the sets $\scrG(3,\Lambda_{N_2})$. 

By adapting the proof of Lemma \ref{L:Deg} with probability at least 
\be
1- \tilde{c} e^{-\log^2 N}
\ee
no vertex has degree exceeding $\log^2 N$ after we sample all edges with length no more than $N_3$. Thus, 
\be
|B_{L(x)}(x)| \leq \log^2 N ||B_{L(x)-1}(x)| \leq N_{4}^{s-d} \log^4 N.
\ee
Let 
\be
\BB^c = \{\exists x \text{ such that } B_{L(x)}(x) \cap \scrG = \varnothing, L(x)< \infty\}
\ee
Then on $\BB$, if $x \in \scrG(3, \Lambda_{N_2})$
\be
d_{\Z^d, \infty} (x, \Lambda_{N_2}) \leq N_3 L(x) \leq N_3 N_4^{s-d}\log^{2} (N)
\ee
thus 
\be
\rho N_2^d \leq \scrG(3, \Lambda_{N_2}) \leq N_2^d + 2d2^{d-1}N_2^{d-1} (N_3N_4^{s-d}\log^{2} (N)).
\ee
Combining the estimates gives the result.
\end{proof}

\subsubsection[Allocation Phase 2]

Now we will work under the conditional probability measure determined by the $\sigma$-algebra $\HH_{N_3}$ and the event $E$ from Lemma \ref{L:P1}.
We must account for the small clusters of $\scrX_{N_3}$ which connect to  $\CC(3, \scrG)$ after revealing the remaining edges between $\CC(3, \scrG)$ and $\scrG^c\cap \scrO$.  Let $ \CC_1, \dotsc, \CC_m$ enumerate the small clusters of $\scrX_{N_3}$.  As was observed above, on $E$ (which is $\HH_{N_3}$ measurable), the small clusters $\CC_1, \dotsc \CC_m$ of $\scrX_{N_3}$ have size at most $N_4^{s-d}\log^4 N$.  Further, the only way to connect these clusters together (or, for that matter, to connect them to $\CC(3, \scrG)$) is via revealing bonds of $\ell^\infty$ length at least $N_3$.  These two conditions prepare the way to dominate cluster sizes by subcritcal branching processes.

Let $\scrX$ denote the graph  
obtained by adding all edges in $\scrO \times (\scrO \cap \scrG^c)$ to the graph $\scrX_{N_3}$.  Let $ \CC_i \overset{\scrX}{\leftrightarrow} \CC(3, \scrG)$ denote the event that $\CC_i$ is connected to $\CC(3, \scrG)$ after all edges outside $\scrG$ have been sampled.  Consider the set $\DD:=\cup_{i: \CC_i \overset{\scrX}{\leftrightarrow} \CC(3, \scrG)}
\CC_i$.  We may define a secondary allocation $\psi: \DD \rightarrow \PP_{N_2}$ using a mechanism entirely analogous to that leading to Definition \eqref{Eq:Part}. Let 
\be
\scrG(0, \Lambda_{N_2})= \scrG(3, \Lambda_{N_2}) \cup \{x \in \DD: \psi(x)= \Lambda_{N_2}\}.
\ee
Note that even without sampling the remaining edges of $\scrG$, $(\scrG(0, \Lambda_{N_2}))_{\Lambda_{N_2} \subset \scrO}$ defines a partition of the largest component of $\scrO$.  Thus it does not depend on edges between $\scrG(\Lambda_{N_2})_{\Lambda_{N_2} \subset \scrO}$ with length larger than $N_2$.

We aim to prove the following:
\begin{lemma}
\label{L:P2}
There exist universal constants $c_7, c_8, \delta_4, \rho_3, \rho_4>0$ so that for all $N> N_1>N_2>N_3> N_4$ chosen as in \eqref{eq:chooseN}, there exists an event $F\subset \Omega$ with 
\be
\mu(F|E, \HH_{N_3}) \geq 1- c_7e^{-c_8 \log^2 N}
\ee on which the following holds:
For all $\Lambda_{N_2} \in  \PP_{N_2}$,
\be
\textrm{Diam}(\scrG(0, \Lambda_{N_2})) \leq N_4^{s-d} \log^{3} N +\textrm{Diam}(\scrG(3, \Lambda_{N_2}))
\ee
and
\be
\rho_3 N_2^d \leq |\scrG(0, \Lambda_{N_2})| \leq \rho_4 N_2^d.
\ee
Also, $\scrG(\Lambda_{N_2}) \subset \scrG(0, \Lambda_{N_2})$ and the $\scrG(0, \Lambda_{N_2})$ form a partition of the largest component of $\scrO$ into connected subsets.
\end{lemma}

\begin{proof}
In this proof we work under the conditional measure $\mu(\cdot|\HH_{N_3}, E)$.  
Denote the set of connected components in $\scrX_{N_3}$ disjoint from the core as $\mathcal S$.  We now reveal all edges not observed in  $\HH_{N_3}$ between clusters in $\mathcal S$. Note that all of these edges must be of length greater than $N_3$.
%\be
%\mathcal S=\{\CC \subset \scrX_{N_3}: \CC \cap \CC(3, \scrG) = \varnothing \}.
%\ee
For each component $\CC \in \mathcal S$ let 
\be
\textrm{deg}({\CC})= \#\{\la x, y \rra: x \in \CC, y \in \CC', \|x-y\|_{\infty} > N_3, \omega_{\la x, y\rra}=1 \text{ for some } \CC' \in \mathcal S \}.
\ee
Then for each $\CC \in \mathcal S$
\be
\label{Eq:Birth}
\E_\mu(\textrm{deg}({\CC})=k| \HH_{N_3}, E) \leq \frac{|\CC|^k}{k!} N_3^{(d-s)k}.
\ee

Consider the course grained graph $CG_{\mathcal S}=(V_{\mathcal S}, \EE_{\mathcal S})$ obtained after revealing all such long edges between the $\CC_j \in \mathcal S$ and then contracting the components of $\mathcal S$ to points.  
Let $\hat{C}(\CC_i)$ denote the connected component of the contraction of $\CC_i$ in this graph.  

Since on $E$ we have a uniform upper bound on $\{|\CC_j|: \CC_j \in \mathcal S\}$ of $N_4^{s-d} \log^2 N$, for each $i$, we can dominate $Diam{\hat{C}(\CC_i)}$ using a subcritical branching process with
mean birth
\be
\eta= N_4^{s-d} \log^{2}(N)N_3^{d-s}
\ee
If $X_n$ denotes a branching process with this mean $\eta$, we couple $X_n$ to $\hat{C}(\CC_i)$ so that births occur whenever we see a new (contraction of) $\CC_k$, exploring $\hat{\CC}(\CC_i)$ starting from $\CC_i$.
As
\be
\BbbP(X_n>0) \leq \eta^n.
\ee
we easily conclude that
\be
\label{Eq:SCS}
\mu(\exists \CC_i: \textrm{Diam}(\hat{C}(\CC_i)) >\log N |E, \HH_{N_3} ) \leq N^d  \eta^{\log N}.
\ee

Again, since the size of each individual cluster $\CC_i$ is bounded by $N_4^{s-d} \log^{2} N$ on $E$, there exist constants $\tilde c_1, \tilde c_2>0$ so that we have, with conditional probability at least $1- \tilde c_1e^{-\tilde c_2\log^2 N}$, 
\be
\label{D:B}
\textrm{Diam}(\cup_{j \in \hat{C}(\CC_i)} \CC_j) \leq N_4^{s-d} \log^{ 3} N
\ee
for all $\CC_i\in \mathcal S$.
Let 
\be
F_1:= \{ \textrm{Diam}(\hat{C}(\CC_i)) <\log N \: \forall \CC_i \in \mathcal S\}.
\ee 
At this point we conclude that no matter how the cluster $\cup_{j \in \hat{C}(\CC_i)} \CC_j$ connects to $\scrG(3, \Lambda_{N_2})$,
 \be
\label{D:B1}
\textrm{Diam}(\scrG(0, \Lambda_{N_2})) \leq N_4^{s-d} \log^{ 3} N + \textrm{Diam}(\scrG(3, \Lambda_{N_2}))
\ee
on $F_1 \cap E$.

Next we consider volume growth.  We can bound $|\hat C(\CC_i)|$ by the total number of offspring $T$ of the process $(X_n)_{n \in\bN}$.  
We have
\be
|\hat{C}(\CC_i)| \leq 1 + T= 1 + \sum_{n=1}^{\infty}X_n
\ee
By \eqref{Eq:Birth} we have that
\be
\E[e^{\lambda X_1}] \leq e^{e^{\lambda} \eta}.
\ee
Thus if $(Y(i))$ denotes an i.i.d. sequence distributed as $X_1$, then
\be
\E[e^{\lambda X_n}]= \E[e^{\lambda \sum_{i=1}^{X_{n-1}} Y(i)}]  \leq \E[e^{ X_{n-1} e^{\lambda}\eta}] 
\ee
Iterating the bound using the natural filtration for $(X_n)_{n=1}^{\infty}$, the monotone convergence theorem implies
\be
\E[e^{\lambda T}] \leq e^{f(\lambda, \eta)}
\ee
where $f(\lambda, \eta)$ is the minimal solution to the equation
\be
f(\lambda, \eta)e^{- f(\lambda, \eta)} = e^{\lambda}\eta
\ee 
provided $\eta$ is sufficiently small (depending on $\lambda$).
Therefore, for any $\lambda>0$, there is $\eta=\eta(\lambda)$ such that if
$N_4^{s-d} \log^{2}(N)N_3^{d-s}< \eta$, then
\be
\E[e^{\lambda |\hat{C}(\CC_i)|}|E, \HH_{N_3}] \leq e^{f(\lambda, \eta)}
\ee
so that we may conclude
\be
\mu( \exists i:  |\hat{C}(\CC_i)| > \log^2N |E, \HH_{N_3}] \leq N^d e^{-\lambda \log^2N} e^{f(\lambda, \eta)}.
\ee
Therefore, we can find constants $\tilde{c_3}, \tilde{c_4}$ so that with conditional probability at least $1-\tilde{c_3}e^{-\tilde{c_4}\log^2 N}$,
\be
Vol(\cup_{j \in \hat{C}(\CC_i)} \CC_j) \leq N_4^{s-d} \log^4 N \quad \forall \CC_i \in \mathcal S.
\ee 

The only thing left to do is make sure we don't allocate too many clusters $\cup_{j \in \hat{C}(\CC_i)} \CC_j$ to the same connected subset $\scrG(3, \Lambda_{N_2})$. 
According to our work in Phase 1 of the allocation, on $E$ we have 
\be
\rho_1 N_2^d \leq |\scrG(3, \Lambda_{N_2})| \leq \rho_2 N_2^d.
\ee  
We have not yet revealed the long edges between $\DD$ and $\scrG(3, \Lambda_{N_2})$ for any $\Lambda_{N_2}$.  This is the next step.
Let
\be
\NN_{\scrG(3, \Lambda_{N_2})}= \sum_{\substack{x \in \scrG(3, \Lambda_{N_2})\\
y \in \scrG^c}} \sum_{\|y-x\|_\infty > N_3} \mvone_{\{\omega_{\la x, y\rra}=1\}}. %%I changed N_{\scrG(3, \Lambda_{N_2})} to \NN_{\scrG(3, \Lambda_{N_2})} to distinguish it from the other N's
\ee
Then $\NN_{\scrG(3, \Lambda_{N_2})}$ is the sum of conditionally independent random variables and in fact is conditionally independent of the connectivity of $\cup_{j\in \hat{C}(\CC_i)} \CC_j$.   If we separate the summands via length, then for each $r \in [N_3, N] \cap \Z$,    we have at most $r^{d-1} \rho_2 N_2^{d}$ independent, approximately identically distributed summands of distribution type $Ber(p= r^{-s})$.
To finish the argument, we use large deviations to bound the number of long edges emanating from any set $\scrG(3, \Lambda_{N_2})$.  
Let 
\be
Var_{N_3}(N_{\scrG(3, \Lambda_{N_2})}) = \E_\mu\left[\left(\NN_{\scrG(3, \Lambda_{N_2})} - \E_\mu[N_{\scrG(3 \Lambda_{N_2})}|E, \HH_{N_3}]\right)^2| E, \HH_{N_3}\right)
\ee denote the conditional variance of $N_{\scrG(3, \Lambda_{N_2})}$.  Then the Azuma-H\"{o}effding inequality implies that for any $\lambda , t> 0$
\begin{multline}
\mu(|N_{\scrG(3, \Lambda_{N_2})} - \E_{\mu}[N_{\scrG(3, \Lambda_{N_2})}|E, \HH_{N_3}]|/\sqrt{Var_{N_3}(N_{\scrG(3, \Lambda_{N_2})})} \geq N_2^{\lambda} | E, \HH_{N_3})\\
 \leq e^{-tN_2^{\lambda}} e^{t^2/2}
\end{multline}
Taking a union bound and optimizing over $t$,
\begin{multline}
\label{Eq:Edge}
\mu(\exists \Lambda_{N_2}: |N_{\scrG(3, \Lambda_{N_2})} - \E_\mu[N_{\scrG(3, \Lambda_{N_2})}|E, \HH_{N_3}]|/\sqrt{Var_{N_3}(N_{\scrG(3, \Lambda_{N_2})})} \geq N_2^{\lambda} | E, \HH_{N_3}) \\
\leq \left(\frac{N}{N_2}\right)^d e^{-N_2^{2\lambda}/2}.
\end{multline}
Since  
\be
\E_\mu[N_{\scrG(3, \Lambda_{N_2})}|E, \HH_{N_3}] \leq C_1 N_2^d N_3^{d-s}
\ee
and
\be
Var_{N_3}(N_{\scrG(3, \Lambda_{N_2})}) \leq C_2 N_2^d N_3^{d-s}
\ee
choosing $N_2^{\lambda}=N_2^{d/2}N_3^{(d-s)/2}$ gives
\be
\label{Eq:Edge1}
\mu(\exists \Lambda_{N_2}: N_{\scrG(3, \Lambda_{N_2})}> 2C_1\vee C_2  N_2^d N_3^{d-s} | E, \HH_{N_3}) \\
\leq \left(\frac{N}{N_2}\right)^d e^{-N_2^{2\lambda}/2}.
\ee

Putting \eqref{Eq:Edge1} together with \eqref{Eq:SCS}, we find
\begin{multline}
\mu(\exists \Lambda_{N_2}: Vol(\scrG(0, \Lambda_{N_2})) \geq \rho_2 N_2^{d} +2C_1 \vee C_2 N_2^d N_3^{d-s} N_4^{s-d}\log^{ 4} N |E, \HH_{N_3})\\
 \leq \left(\frac{N}{N_2}\right)^d e^{-N_2^dN_3^{d-s}} + N^d e^{- \lambda \log^2 N} e^{f(\lambda, \eta)}
\end{multline}
\end{proof}

Before moving to the the derivation of heat kernel bounds, let us pause to record the last statement of Theorem \ref{t:diameter} as a corollary to the previous Lemma:
\begin{corollary}
\label{C:connect}
If $\CC^1(N)$ denotes the largest component in $B_N(0):=[-N, N]^d$, there is $\epsilon>0$ so that
\be
 \mu \left( 0 \leftrightarrow  B^c_{N}(0)| 0\notin \CC^1(N) \right)
 \leq  C N^{- \epsilon}.
\ee
\end{corollary} 
\begin{proof}
This is a simple application of the branching process argument exploited in the previous lemma.  The proof is omitted.
\end{proof}

\subsection{Step 3: Assembling the Estimates}
Next we derive the required spectral gap bounds.  Let $\HH_+$ denote the $\sigma$-algebra generated by the construction so far, that is $\FF_{N_2}$, plus all edges between vertices in $\scrG^c\cap \scrO$ - leaving out only edges inside the core of length greater than $N_2$.
Let us sum up our construction to this point:
There exist universal constants $c_9, c_{10}, \delta_5, \rho_5, \rho_6> 0$ so that if
$(N_i)_{i=0}^4$ are chosen as \eqref{eq:chooseN} then we can find an event $R$ with $\mu(R) \geq 1- c_9 e^{-c_{10} \log^2 N}$, measurable with respect to $\HH_+$ and so that on  $R$:
\begin{enumerate}
\item
\label{II:1}
$\scrG$ exists and satisfies the properties of Lemma \ref{L:CoreProps}.
\item
\label{II:2}
The connected components of $\scrO$ are measurable with respect to $\HH_+$.  In other words, adding the states of edges which have not yet been revealed does not change the connectivity.  
\item
\label{II:3}
On $R$ the small components of $\scrO$ have size at most $N_4^{s-d}\log^4N$.  
\item
\label{II:4}
The largest component $\scrM$ of $\scrO$ may be partitioned into connected random subsets $(\scrG(0, \Lambda_{N_2}))_{\Lambda_{N_2} \in \PP_{N_2}}$ so that
\begin{itemize}
\item $\textrm{Diam}(\scrG(0, \Lambda_{N_2}) \leq \log^{\delta_5} N$
\item $\rho_5 N_2^d \leq |\scrG(0, \Lambda_{N_2})| \leq \rho_6 N_2^d$.
\item $\scrG(\Lambda_{N_2}) \subset \scrG(0, \Lambda_{N_2})$
\end{itemize}
\item
\label{II:5} 
Finally edges of length at least $N_2$ \textit{between} vertices in $\scrG$ have not been sampled and are thus conditionally independent of $R$.
\end{enumerate}

At this point we are in a position to apply the strategy of \cite{BBYr} to obtain the bounds on gaps of individual sets associated with the $ \Lambda_{N_1}$ which tile $\scrO$. 
The idea is as follows:  for each $\Lambda_{N_1}$, consider $\{\scrG(0, \Lambda_{N_2}): \Lambda_{N_2} \subset \Lambda_{N_1}\}$.  We now sample the edges between the corresponding $\scrG(0,\Lambda_{N_2})$.  
Define
\be
\scrG(0, \Lambda_{N_1}):= \cup_{\Lambda_{N_2} \subset \Lambda_{N_1}} \scrG(0, \Lambda_{N_2})
\ee
along with these newly sampled edges.
We will bound the spectral gap of $\scrG(0, \Lambda_{N_1})$ using multi-commodity flows \cite{Sinclair}.

\textbf{A Brief Primer on Multicommodity Flows (exposition taken from \cite{BBYr}):}
Let $P$ be the transition matrix of a reversible Markov chain,
with stationary distribution $\pi$.
Let $V$ be the set of states of the chain, and let $\EE$ be the set of
oriented edges; i.e.
$$\EE = \set{ (x,y) \in V \times V \ : \ P(x,y) > 0 } . $$

For $x,y \in V$ let $\Gamma(x,y)$ be the set of all simple paths from
$x$ to $y$. Let $\Gamma = \cup_{x \neq y \in V} \Gamma(x,y)$.

A \emph{flow} is a function
$f:\Gamma \to [0,1]$ such that for all $x,y \in V$
$$ \sum_{\gamma \in \Gamma(x,y)} f(\gamma) = \pi(x) \pi(y) . $$
The \emph{edge load} of an edge $e \in E$ is defined as
$$f(e) = \sum_{\substack{ \gamma \in \Gamma \\ \gamma \ni e}} f( \gamma ) | \gamma | $$
where $|\gamma|$ denotes the number of edges in $\gamma$.
The \emph{congestion} of a flow $f$ is defined as
$$ \rho(f) = \max_{(a,b) \in \EE} \frac{1}{\pi(a) P(a,b)} f((a,b)) . $$

Theorem 5' of \cite{Sinclair} states that if the eigenvalues of
$P$ are $1 > \lambda \geq \lambda_3 \geq \cdots \geq \lambda_n$
(where $n=|V|$), then for any flow $f$, $(1- \lambda)^{-1} \leq \rho(f)$.  Furthermore, Theorem 8 in \cite{Sinclair} shows that if $P$
induces an ergodic Markov chain (i.e. if $\lambda_n > -1$), then
there exists a flow $f^*$ such that
$\rho(f^*) \leq 16 \tau$, where $\tau$ is the mixing time of the chain.
We call $f^*$ the \emph{optimal flow} for $P$.

The following (among other things) was proved in \cite{Eyal}:
\begin{theorem}[Benjamini et. al. Theorem 1.2]
\label{BET}
There exists $C>0$
such that if $G$ is chosen according to $G(n, p)$ with $p\geq C\log n/ n$ then
the mixing time of $G(n, p)$ (which in this regime is connected) has a bound $\tau_{G(n, p)} = O(\log n)$ with probability tending to $1$ as $N \ra \infty$.
\end{theorem}
 %%%%%%%%%%%%Do we need the isoperimetric constant

For us, $G$ will be one of the graphs $\scrG(0, \Lambda_{N_1})$.  If
\be
\textrm{deg}_{\scrG(0, \Lambda_{N_1})}(x)= |\{y \in \scrG(0, \Lambda_{N_1}): \omega_{\la x,y \rra}=1\} |
\ee
Then we let
\be
\pi_{\scrG(0, \Lambda_{N_1})}(x)= \frac{\textrm{deg}_{\scrG(0, \Lambda_{N_1})}(x)}{\sum_{y \in \scrG(0, \Lambda_{N_1})}\textrm{deg}_{\scrG(0, \Lambda_{N_1})}(y)}.
\ee
Of course, this is the stationary measure for SRW on $\scrG(0, \Lambda_{N_1})$.  Let $\textrm{Gap}_{\scrG(0, \Lambda_{N_1})}$
denote the spectral gap associated to this Markov chain.

A word about notation below.  We will consider various graphs induced by long range percolation and random graph processes coupled to these induced graphs.  Given such a graph $G$, we will use $V(G)$ to refer to the vertices of $G$ and $\EE(G)$ will denote the set of edges of G.
\begin{lemma}[Lower Bound on Spectral Gap for $\scrG(0, \Lambda_{N_1})$]\label{l:gap}
There exist universal constants $c_{11}, c_{12}, \delta_6>0$ so that if $N> N_1> N_2>N_3>N_4$ be fixed as in \eqref{eq:chooseN} then
\be
\mu(\textrm{Gap}_{\scrG(0, \Lambda_{N_1})}  >   c N^{d-s}_1/\log^{\delta_6} N 
 \textrm{ for all } \Lambda_{N_1} \in \PP
|R, \HH_+  ) \geq 1 - c_{11}\left(N/N_1\right)^d e^{-c_{12} \log^2 N_1}  
\ee
\end{lemma}
\begin{proof}
Due to the conditional independence built into our construction, if $d_{\Z^d, \infty}(\Lambda_{N_2}, \Lambda_{N_2}') \leq k N_2$ then for each $x \in \scrG(\Lambda_{N_2})$
\be
\label{Eq:Key1}
\mu(x \nleftrightarrow \scrG(\Lambda_{N_{2}}')|R, \HH_+ ) \leq e^{-\rho^2 (k+1)^{-s} N_2^{d-s}}.
\ee
By construction, $|\scrG(\Lambda_{N_2})|\geq \rho_2 N_1^{s-d}\log^{\gamma_1} N_1$ and so
\be
|\scrG(\Lambda_{N_2})| (k+1)^{-s} N_2^{d-s} \geq (N_2/N_1)^{d}= N_1^{s-2d} \log^{3} N_1
\ee
Fix a block $\Lambda_{N_1}\in \PP_{N_1}$. %%%Will now suppress the \Lambda_{N_1} from the notation in this lemmas
We may thus compare the block connectivity of $(\scrG(\Lambda_{N_2}))_{\Lambda_{N_2} \subset \Lambda_{N_1}}$ to an Erd\"{o}s-R\'{e}nyi graph $G(n, p)$ of size $n=N_1^{2d-s}/\log^{3} N_1$ with $p=N_1^{s-2d} \log^{3} N_1 =(\log^6 N_1)/n$.  The comparison graph $G(n, p)$ is, by our choice of $N_2$, well in the super critical range.

To increase the probability of our graph having the required properties to $1-e^{-c\log^2 N}$ we use the same amplification technique that was employed in Lemma \ref{L:Diam}.  We may view the course grained block percolation process as containing $\lfloor \log^{2} N_1/C\rfloor$ independent identically distributed samples each distributed as $G(n, p')$ (where $C$ has been chosen sufficiently large so that Theorem \ref{BET} holds)
and $p'=C N_1^{s-2d} \log N_1 = C(\log^4 N_1)/n$. Let $(ER_j)_{j=1}^{\lfloor \log^{2} N_1/C\rfloor}$ denote the course-grained i.i.d. copies and let $\scrE \scrR_j(\Lambda_{N_1})$ denote the Erd\"{o}s R\'{e}nyi random graph samples. 
From Theorem \ref{BET}, we immediately conclude that except with probability at most $e^{-\lfloor \log^{2} N_1/C\rfloor}$, at least one of the comparison Erd\"{o}s-R\'{e}nyi graphs has a mixing time of order
\be
\label{ER-bound}
\tau \leq c_{12} \log N.
\ee

%\textbf{Check and adjust if nesc:}
%Let 
%\be
%\scrG_j(0, \Lambda_{N_1})= \cup_{\Lambda_{N_2} \subset \Lambda_{N_1}} \scrG(0, \Lambda_{N_2}) \cup ER_j(\Lambda_{N_1})
%\ee
%By the theory of Dirichlet forms, the spectral gap of SRW on $\scrG(0, \Lambda_{N_1})$ satisfies
%\be
%\lambda_{\scrG(0, \Lambda_{N_1})} \geq \max_i \lambda_{\scrG_i(0, \Lambda_{N_1})}.
%\ee
%and we may obtain a lower bound by analyzing $\scrG_j(0, \Lambda_{N_1})$ for any $j$.

We choose any $j_0$ which satisfies \eqref{ER-bound} and define a multicommodity flow on $\scrG_{j_0}(0, \Lambda_{N_1})$ 
%%% A: do we need the subscript j here?
%%N:  Strictly speaking prob. not, but the flow we construct is concentrated on this graph
as in \cite{BBYr}, also described in detail below.
Thanks to Sinclair's work \cite{Sinclair}, for $\scrE \scrR_{j_0}$ the coupled Erd\"{o}s-R\'{e}nyi graph has on optimal flow $f_{ER}$ with
$\rho(f_{ER}) \leq 16c \log N_1$.

The lower bound on the spectral gap of $\scrG(0, \Lambda_{N_1})$, 
is derived by constructing a flow supported on $\scrG_{j_0}(0, \Lambda_{N_1})$ using the optimal flow $f_{ER}$.  We follow \cite{BBYr} rather closely.

Let $\pi_{\scrG(0, \Lambda_{N_1})},\pi_{\scrE \scrR_{j_0}}$ denote the stationary distribution of 
SRW on the respective graphs $\scrG(0, \Lambda_{N_1}), \scrE \scrR_{j_0}(\Lambda_{N_1})$.
Let 
\be
\Gamma(\scrG_{j_0}(0, \Lambda_{N_1})),\Gamma(\scrE \scrR_{j_0}(\Lambda_{N_1}))
\ee 
be the set of simple paths in $\scrG_{j_0}(0, \Lambda_{N_1}),\scrE \scrR_{j_0}(\Lambda_{N_1})$,
and let
\be
\Gamma(x,y;\scrG_{j_0}(0, \Lambda_{N_1})), \Gamma(i,k;\scrE \scrR_{j_0}(\Lambda_{N_1}))
\ee 
be the set of simple paths in 
\be
\scrG_{j_0}(0, \Lambda_{N_1}),\scrE \scrR_{j_0}(\Lambda_{N_1})
\ee 
from $x$ to $y$ and $i$ to $k$, respectively.
For a path $\gamma$, let $\gamma^+$ be the starting vertex of $\gamma$, and let $\gamma^-$ be the ending 
vertex of $\gamma$ (specifically for edges $e = (e^+,e^-)$).

For each $i \in \scrE \scrR_{j_0}$, let ${\Lambda_{N_2}(i)}$ denote the block in $\scrG_{j_0}(0, \Lambda_{N_1})$ coupled to $i$.
For $(i,k) \in \EE(\scrE \scrR_{j_0})$ let $e(i,k)$ be a specific edge of $\EE(\scrG_{j_0}(0, \Lambda_{N_1}))$ given by the coupling with $x \in \scrG(\Lambda_{N_2}(i))$ and $y \in \scrG(\Lambda_{N_2}(k))$
(by definition, under our coupling there always exists at least one such edge).

For each pair $x,y \in \scrG(0, \Lambda_{N_2})$ let $\gamma(x,y)$ be a path in $\scrG(0, \Lambda_{N_2})$ that realizes the graph
distance between $x$ and $y$ in $\scrG(0, \Lambda_{N_2})$ (i.e. a geodesic).  In case $x=y$ let $\gamma(x,x)$
be the empty path.

For $\eta \in \Gamma(i,j;\scrE \scrR_{j_0}(\Lambda_{N_1}))$, and $x \in V(\scrG(0, \Lambda_{N_2}(i))), y \in V(\scrG(0, \Lambda_{N_2}(k)))$, 
define 
$\gamma(\eta,x,y) \in \Gamma(x,y;\scrG_{j_0}(0, \Lambda_{N_1}))$ by interpolating $\eta$ using the specified edges $e(i,k)$
and geodesics $\gamma(x,y)$.  In other words, if $\eta = e_1 e_2 \cdots e_{|\eta|}$,
then
$$ \gamma(\eta,x,y) = \gamma(x,e_1^+) e(e_1^+,e_1^-) \gamma(e_1^-,e_2^+) e(e_2^+,e_2^-)
\cdots e(e_{|\eta|}^+,e_{|\eta|}^-)  \gamma(e_{|\eta|}^-,y) . $$
Setting 
$\Delta = \max_{\Lambda_{N_2} \subset \Lambda_{N_1}}\mathrm{Diam}(\scrG(0, \Lambda_{N_2}))$
we get that $|\gamma(\eta,x,y)| \leq (\Delta+1) |\eta|$.

As mentioned above,
by Theorem 8 of \cite{Sinclair},
there exists a constant $c_1>0$ such that
\begin{align} \label{eqn:flow Gamma'}
\forall & \ (i,k) \in \EE(\scrE \scrR_{j_0}(\Lambda_{N_1})) \qquad\\
&|\EE(\scrE \scrR_{j_0}(\Lambda_{N_1}))| \sum_{\substack{ \eta \in \Gamma(\scrE \scrR_{j_0}(\Lambda_{N_1})) \\ \eta \ni (i,k) } }
f_{ER}(\eta) |\eta|  \leq 16 \tau(\scrE \scrR_{j_0}(\Lambda_{N_1})) 
\leq c_1 \log N . \nonumber
\end{align}

We now define the flow $f$ on $\scrG(0, \Lambda_{N_1})$.
Let $x,y \in V(\scrG(0, \Lambda_{N_1}))$, and let $i,k$ be such that $x \in \scrG(0, \Lambda_{N_2}(i))$ and $y \in \scrG(0, \Lambda_{N_2}(k))$.

If $i=k$ route all the flow along $\gamma(x,y)$ so $f(\gamma(x,y)) = \pi_{\scrG(0, \Lambda_{N_1})}(x) \pi_{\scrG(0, \Lambda_{N_1})}(y)$.

On the other hand, if $i \neq k$, then for any $\eta \in \Gamma(i,k;\scrE \scrR_{j_0}(\Lambda_{N_1}))$, $x \in \scrG(0, \Lambda_{N_2}(i))$ and $y \in \scrG(0, \Lambda_{N_2}(k))$ set 
$$ f(\gamma(\eta,x,y)) = \frac{f_{ER}(\eta)}{\pi_{\scrE \scrR_{j_0}(\Lambda_{N_1})}(i) \pi_{\scrE \scrR_{j_0}(\Lambda_{N_1})}(k)} 
\cdot \pi_{\scrG(0, \Lambda_{N_1})}(x) \pi_{\scrG(0, \Lambda_{N_1})}(y) , $$
and $0$ otherwise.

We bound the congestion of $f$  along some edge $(x,y) \in \EE(\scrG_{j_0}(0, \Lambda_{N_1}))$. Let $i,k$ be such that $x \in \scrG(0, \Lambda_{N_2}(i))$ and $y \in \scrG(0, \Lambda_{N_2}(k))$.

{\bf Case 1:} $i \neq k$.
In this case, any path $\gamma$ that contains the edge $(x,y)$,
such that $f(\gamma) > 0$ 
must be of the form $\gamma = \gamma(\eta,z,w)$ for some $\eta \in \gamma(\scrE \scrR_{j_0}(\Lambda_{N_1}))$ that contains $(i,k)$.
Thus, 
\begin{align} 
&\sum_{\substack{\gamma \in \Gamma(\scrG_{j_0}(0, \Lambda_{N_1})) \\ \gamma \ni (x,y)} } f(\gamma) |\gamma| \nonumber\\
& \quad \leq 
\sum_{\substack{\eta \in \Gamma(\scrE \scrR_{j_0}(\Lambda_{N_1})) \\ \eta \ni (i,k) } }
\sum_{\substack{ z \in \scrG(0, \Lambda_{N_2}(\eta^+)) \\ w \in \scrG(0, \Lambda_{N_2}(\eta^-))} } 
\frac{\pi_{\scrG(0, \Lambda_{N_1})}(z) \pi_{\scrG(0, \Lambda_{N_1})}(w)}{\pi_{\scrE \scrR_{j_0}(\Lambda_{N_1})}(\eta^+) \pi_{\scrE \scrR_{j_0}(\Lambda_{N_1})}(\eta^-) } \cdot 
f_{ER}(\eta) |\eta| (\Delta+1). \nonumber
\end{align}
It easily follows that the right hand side is bounded by
\begin{equation}
\label{Eq:Int}
(\Delta+1) \cdot (\max_\ell \pi_{\scrG(0, \Lambda_{N_1})}(\scrG(0, \Lambda_{N_2}(\ell))^2) \cdot |\EE(\scrE \scrR_{j_0}(\Lambda_{N_1}))|^2 \cdot 
\sum_{\substack{\eta \in \Gamma(\scrE \scrR_{j_0}(\Lambda_{N_1})) \\ \eta \ni (i,j) } } f_{ER}(\eta) |\eta|.
\end{equation}
Next, using \eqref{eqn:flow Gamma'},
we can bound \eqref{Eq:Int} by
\begin{equation}
\label{eqn:long edge}
c_1 \log (N) \cdot (\Delta+1) \cdot (\max_\ell \pi_{\scrG(0, \Lambda_{N_1})}(\scrG(0, \Lambda_{N_2}(\ell))^2) \cdot |\EE(\scrE \scrR_{j_0}(\Lambda_{N_1}))|.
\end{equation}

{\bf Case 2:} $i=k$.
In this case, any path $\gamma$ that contains the edge $(x,y)$,
such that $f(\gamma) > 0$, 
is one of the follwing:
Either it is of the form $\gamma = \gamma(\eta,z,w)$ for some $\eta \in \Gamma(\scrE \scrR_{j_0}(\Lambda_{N_1}))$ that
contains the vertex $i$, or it is of the form $\gamma = \gamma(z,w)$ for some $z,w \in \scrG( 0,\ \Lambda_{N_2}(i))$.
Any path $\eta \in \Gamma(\scrE \scrR_{j_0}(\Lambda_{N_1}))$ that contains the vertex $i$ must contain some edge
$(i,k) \in \EE(\scrE \scrR_{j_0}(\Lambda_{N_1}))$.  Thus, using \eqref{L:Deg}
and \eqref{eqn:long edge},
\begin{align} \label{eqn:short edge}
&\sum_{\substack{\gamma \in \Gamma(\scrG_{j_0}(0, \Lambda_{N_1})) \\ \gamma \ni (x,y)} } f(\gamma) |\gamma| \nonumber \\
&\quad  \leq \sum_{z,w \in \scrG(\Lambda_{N_2}(i))} f(\gamma(z,w)) |\gamma(z,w)| \nonumber\\
&\quad + \sum_{k:(i,k) \in \EE(\scrE \scrR_{j_0}(\Lambda_{N_1}))}
\sum_{\substack{\eta \in \Gamma(\scrE \scrR_{j_0}(\Lambda_{N_1})) \\ \eta \ni (i,k) } } 
\sum_{\substack{ z \in \scrG(\Lambda_{N_2}(\eta^+)) \\ w \in \scrG(\Lambda_{N_2}(\eta^-)) } }
f(\gamma(\eta,z,w)) |\gamma(\eta,z,w)|
\nonumber \\
& \quad \leq 
\Delta \sum_{z,w \in \scrG(\Lambda_{N_2}(i))} \pi_{\scrG(0, \Lambda_{N_1})}(z) \pi_{\scrG(0, \Lambda_{N_1})}(w) \nonumber
\\
&\quad + \sum_{k:(i,k) \in \EE(\scrE \scrR_{j_0}(\Lambda_{N_1}))}  
\sum_{\substack{\eta \in \Gamma(\scrE \scrR_{j_0}(\Lambda_{N_1})) \\ \eta \ni (i,k) } } 
\frac{\pi_{\scrG(0, \Lambda_{N_1})}(\scrG(0,\Lambda_{N_2}(\eta^+))) \pi_{\scrG(0, \Lambda_{N_1})}(\scrG(0, \Lambda_{N_2}(\eta^-)))}{\pi_{\scrE \scrR_{j_0}(\Lambda_{N_1})}(\eta^+) \pi_{\scrE \scrR_{j_0}(\Lambda_{N_1})}(\eta^-) } (\Delta+1)
f_{ER}(\eta) |\eta|
\nonumber \\
& \quad  \leq \Delta \pi_{\scrG(0, \Lambda_{N_1})}( \scrG(0, \Lambda_{N_2}(i)))^2  \nonumber\\
& \quad  \quad + c_1 \log (N) 
\cdot  \max_i \deg_{\scrE \scrR_{j_0}(\Lambda_{N_1})}(i) 
\cdot  (\Delta+1) \cdot (\max_\ell \pi_{\scrG(0, \Lambda_{N_1})}(\scrG(0, \Lambda_{N_2}(\ell))^2)
\cdot |\EE(\scrE \scrR_{j_0}(\Lambda_{N_1}))| .  %%%%  Feel free to change it back but I found this equation hard to follow
\end{align}

The following all hold with probability at least $1-c'e^{-c'\log^2 N }$:
\begin{itemize}
\item By construction, we have
$\textrm{Diam}(\scrG(0, \Lambda_{N_2})) \leq \log^{\delta_5} N$ for all $\Lambda_{N_2}\in \PP_{N_2}$.
\item By Lemma \ref{L:Deg}, $Deg(\scrO) \leq \log^2 N$.
\item By Lemma \ref{L:P2}
\[
\max_\ell \pi_{\scrG(0, \Lambda_{N_1})}(\scrG(0, \Lambda_{N_2}(\ell)) = \frac{\max_\ell |\EE(\scrG(0, \Lambda_{N_2}(\ell))|}{|\EE (\scrG(0, \Lambda_{N_1}))|}\leq 
\frac{\rho_4N_2^d Deg(\Lambda_1)}{\frac12 \rho_2 N_1^d},
\]
and
\[
|\EE (\scrG(0, \Lambda_{N_1}))| \geq \frac12 |\EE (\scrG(0, \Lambda_{N_1}))|  \geq \frac12 \rho_2 N_1^d
\]
\item Finally by standard concentration results
$Deg(\scrE \scrR_{j_0}) \leq 2C \log^4 N$ and
\[
 |\EE(\scrE \scrR_{j_0}(\Lambda_{N_1}))|  \leq 2n^2 p' = 2N_1^{2d-s}.
\]
\end{itemize}
Putting these estimates together, we have  with probability at least $1-c'e^{-c'\log^2 N }$
\be
\rho(f) \leq  N_1^{s-d} \log^{\delta_6} N
\ee
Finally, applying Theorem 5' of \cite{Sinclair} to the flow $f$ completes the lemma.
\end{proof}

\begin{proof}[Theorem \ref{T:gap}]
The volume bounds are established by Lemmas \ref{L:CoreProps} and \ref{L:P2}.  The bound on the diameter follows by the facts that on the event $R$ each $\Lambda_{N_2} \in \PP_{N_2}$ has diameter at most $\log^{\delta_5} N$ and by the proof of Lemma \ref{L:Diam}.  Finally Lemma \ref{l:gap} establishes the bound on the spectral gap.
\end{proof}

This theorem and the accompanying construction gives us a method of partitioning the largest component of $\scrO$ into connected sets which are approximately blocks but respect the graphs stucture and have good volume, diameter and spectral gap bounds.  Moreoever, taking $N=N_1$ establishes Theorems \ref{t:diameter} and \ref{t:gap_bounds}.

\section{An Abstract Continuous Time Heat Kernel Estimate}
\label{S:HK}
In this section we prove heat kernel upper bounds on the basis of certain \textit{ a priori} hyptheses.   At this point we leave the setting which motivates this paper, LRP.  The goal is to formulate conditions on a graph $G=(V, \EE)$ under which we may prove on-diagonal heat kernel upper bounds via the technology that appears in \cite{Barlow-HKP}, {see also the references there in}.  For convenience, we will work only with the continuous time SRW.  Standard coupling arguments may be employed to obtain discrete time heat kernel bounds as well (see, for example, the appendix to \cite{Berger-Biskup}).
 
Let us fix a (possibly infinite) connected graph $\GG=(V, \EE)$.   Let $\mu(x)=\textrm{\textrm{deg}(x)}$ denote the (possibly non-normalizable) measure on $\GG$ which weights each vertex by its degree.  Let $\LL$ denote the generator of the continuous time simple random walk on $\GG$, normalized to have unit holding times:
\be
\LL f(y) = \sum_{z\in G} \frac{\mvone_{\la y, z\rra \in\EE(\GG)}}{\textrm{deg}(y)} (f(z)-f(y))
\ee
Then  $\LL$ is a self adjoint operator with respect to the Hilbert space $L^2(V, \mu)$
and
\be
(f, \LL f)_{\mu} =- \frac{1}{2} \sum_{y, z \in V} \mvone_{\la y, z\rra \in\EE(\GG)} (f(z)-f(y))^2
\ee
Let $P_t(x, y)$ denote the corresponding transition kernel.

For any \textit{finite} connected subgraph $\HH \subset \GG$ and any vertex $x \in \HH$ let $\textrm{deg}_{\HH}(x)$ denote the degree of $x$ within $\HH$.  Let
\be
\nu_{\HH}(x)=\textrm{deg}_{\HH}(x)/\ZZ(\HH)
\ee
with $\ZZ(\HH)=\sum_{y \in \HH} \textrm{deg}_{\HH}(y)$.

Assume that for each $s \in [T_1, T_2]$, there exists a distinguished connected subset $B(s)$ and a partition $\PP_s$ of $B(s)$ into connected sets $\{\HH: \HH \in \PP_s\}$ and a pair of postive functions $\lambda_s, V_s$ on $\R^+$ (we assume $\lambda_s$ decreases and $V_s$ increases) and a family of universal constants $\{c_i, C_i\}_{i=1}^4$ so that the following assumptions hold:
\begin{enumerate}
\item
\label{I:1}
For all $\HH \in \PP_s$, the spectral gap $\textrm{Gap}_{\HH}$ of the restriction of $\LL$ to $\HH$ satisfies
\be
\textrm{Gap}_{\HH} \geq \lambda_s.
\ee
\item
\label{I:2}
For all $\HH \in \PP_s$, the volume of $\HH$ satisfies
\be
c_1 V_s \leq Vol(\HH) \leq C_1 V_s.
\ee
\item
\label{I:3}
Suppose that for some $\gamma, \tilde{\delta}_1 >0$, the function $\lambda_s$ is linked to $V_s$ by the condition
\be
c_2 V^{-\gamma}_s \log^{-\tilde{\delta}_1} V_s \leq \lambda_s.
\ee
\item
\label{I:4}
Further, let $\Delta_{\PP_s} = \min_{x \in B(s)} \frac{\textrm{deg}_{\HH}(x)}{\textrm{deg}(x)}$
and suppose that there exist constants $c_3, {\tilde{\delta}_2}\ > 0$ so that
\be
\Delta_{\PP_s} \geq c_3 \log^{-\tilde{\delta}_2} V_s.
\ee
\item
\label{I:5}
Next let us suppose that there exists a subset $B_R \subset B(s)$ so that
\be
\sup_{x \in B_R}P_s(x, B(s)^c) \leq C_3 \frac{\log V_s}{V_s }.
\ee
\item
\label{I:6}
Suppose that there exist $ C_4> 0$ so that
\be
2+C_3 \leq  \psi_s V_s \log^{-1} V_s \leq C_4.
\ee
\end{enumerate}

Then we have the following abstract heat kernel bound:
\begin{lemma}
\label{L:AHK}
Under the assumptions (\ref{I:1})--(\ref{I:6}), consider
\be
\psi_t = P_{2t}(x,x)/\textrm{deg}(x)
\ee
for $t \in [T_1/2, T_2/2]$
and $x \in B_R$.

Let $\delta= 2+\tilde\delta_1+ \tilde\delta_2 +\gamma$.
Then there exist $C_5, C_6>0$ (depending only on $\gamma$ and the constants $c_1, C_1, c_2, c_3, C_3, c_4, C_4$),
\be
\psi_t(x) \leq \psi_{T_1}(x)\wedge C_5 \frac{1+ C_6(t-T_1/2))^{-1/\gamma}}{|\log(1+  C_6(t-T_1/2))|^{\delta/\gamma}}
\ee
{when}\ $t \in[ T_1/2, T_2/2]$.
\end{lemma}
\begin{proof}
Our various hypotheses will be explained over the course of the proof.
Following Barlow, let $f_t(y)= \frac{P_t(x, y)}{\mu(y)}$.
Then
\be
\psi_t= \sum_{y \in \GG} f^2_t(y) \mu(y)=(f_t, f_t)_{\mu}.
\ee
Differentiating,
\be
\partial_s \psi_s = -  \sum_{y, z \in V} \mvone_{\la y, z\rra \in\EE(G)} (f_s(z)-f_s(y))^2.
\ee
It is convenient to work with $-\partial_s \psi_s$.  By assumption,
\be
-\partial_s \psi_s \geq \sum_{\HH \in \PP_s}  \sum_{y, z \in \HH} \mvone_{\la y, z\rra \in\EE(G)} (f_s(z)-f_s(y))^2.
\ee

Since
\be
\sum_{y, z \in \HH} \mvone_{\la y, z \rra \in \GG} (f_s(y)-f_s(z))^2
\ee
is the Dirichlet form of $f_s$ for SRW on $\HH$, assumption (\ref{I:1}) implies
\be
\sum_{y, z \in \HH} \mvone_{\la y,z \rra \in \GG} (f_s(y)-f_s(z))^2 \geq  \lambda _s \textrm{Var}_{\HH}(f_s)
\ee
where
\be
\textrm{Var}_{\HH}(f_s):= \sum_{y \in\HH} \textrm{deg}_{\HH}(y) (f_{s}(y) - \E_{\HH}(f_s))^2.
\ee
and
\be
\E_{\HH}(f_s)= \sum_{y\in \HH} \nu_{\HH}(y) f_s(y).
\ee
Thus
\be
-\partial_s \psi_s \geq \lambda_s \sum_{\HH \in \PP_s}  \sum_{y \in\HH} \textrm{deg}_{\HH}(y) f^2_{s}(y) -\frac{1}{\ZZ(\HH)} (\sum_{y\in \HH}\textrm{deg}_{\HH}(y) f_s)^2.
\ee
Since $\sum_{y} \textrm{deg}(y)f_{s}(y) \leq 1$,
\be
\sum_{\HH \in \PP_s} \ZZ(\HH)^{-1}(\sum_{y \in \HH} \textrm{deg}_{\HH}(y) f_s(y))^2 \leq \max_{\HH \in \PP} \ZZ(\HH)^{-1}
\ee
and
\begin{multline}
 \lambda_s \sum_{\HH \in \PP_s}  \sum_{y \in\HH} \textrm{deg}_{\HH}(y) f^2_{s}(y) -\frac{1}{\ZZ(\HH)} (\sum_{y\in \HH}\textrm{deg}_{\HH}(y) f_s)^2\\
 \geq  \lambda_s \sum_{\HH \in \PP_s}  \sum_{y \in\HH} \textrm{deg}_{\HH}(y) f^2_{s}(y) - \lambda_s\max_{\HH \in \PP_s} \ZZ(\HH)^{-1}\\
 \geq
\Delta_{\PP_s} \lambda_s ( \sum_{y \in B(s)} \textrm{deg}(y) f^2_{s}(y) -\Delta_{\PP_s}^{-1} \max_{\HH \in \PP_s} \ZZ(\HH)^{-1})
\\
=\Delta_{\PP_s} \lambda_s ( \psi_s -\Delta_{\PP_s}^{-1} \max_{\HH \in \PP_s} \ZZ(\HH)^{-1})- \ \Delta_{\PP_s} \lambda_s \left( \sum_{y \in B(s)^c} \textrm{deg}(y) f^2_{s}(y)\right)
\end{multline}

Using reversibility and the fact that probabilities are bounded by $1$, assumptions \eqref{I:2} and \eqref{I:5} imply that
\be
\label{Eq:HK}
-\partial_s \psi_s \geq  \Delta_{\PP_s}\lambda_s\left(\psi_s- \frac{(1+C_3)\log V_s}{V_s}\right)
\ee
This explains some of our hypotheses:   by assumption \eqref{I:6}, $\psi_s V_s/\log V_s \geq C_3 +2$.
According to hypotheses \eqref{I:3} and \eqref{I:4}
\be
\label{Eq:HK1}
-\partial_s \psi_s \geq c' \psi_s  V_s^{-\gamma}/\log^{2+\delta_1+ \delta_2} V_s.
\ee
Under the further hypothesis $\psi_s V_s/\log V_s \leq C_4$, we thus obtain
\be
\label{Eq:HK2}
-\partial_s \psi_s \geq c'' \psi_s^{1+\gamma}/|\log \psi_s|^{\delta}
\ee
where $\delta= 2+\delta_1+ \delta_2 +\gamma$.

Using the change of variables $u_s=\psi_s^{-\gamma}$:
\be
\label{Eq:HK3}
\partial_s u_s \geq c'''(\gamma)/\log^{\delta}
 u_s
\ee
where we note $u_s \geq 1$ (and increasing) and $c'''(\gamma)> 0$ is a constant depending only on the exponent $\gamma$.
Integration by parts gives
\begin{multline}
\label{Eq:HK4}
u_s  \log^{\delta}
 u_s \geq c'''(\gamma) s + u_0 \log^{\delta} u_0  + c'''(\gamma) \delta \int_0^s  \log^{\delta-1} u_s\\ \geq c'''(\gamma) s +  u_0
 \log^{\delta} u_0\end{multline}
since $u_s \geq 1 \:  \forall s$.

As $f(u)=   u\log^{\delta}
 u$ increases for $u\geq 1$, solving the equation
\be
 u \log^{\delta}
 u = c''' s
\ee
 for $u$ gives us a lower bound on $u_s$.
 Letting $g(s) = c'''s/ \log^{\delta}(c'''s)$ gives
 \be
 u_s \geq g(s)
 \ee
 by direct calculation.
 
Recalling that $u_s= \psi^{-\gamma}_s$, the claim follows by appropriate choice of the constant $C_5$ since $\psi_s$ is uniformly bounded by $1$.
\end{proof}

\section{Proof of Theorem \ref{T:HKLRP}}

To apply Lemma \ref{L:AHK} we need to check that the assumptions of that lemma hold.   Assumptions \eqref{I:1}-\eqref{I:4} will all follow by Theorem~\ref{T:gap}.  Assumption $\eqref{I:6}$ is a more of a calibration condition for our choice of scale $\PP_s$ than a stringent requirement.
Thus the only assumption that needs further proof is \eqref{I:5}, which will be given next.  Finally, in the subsequent subsection, we will gather all these results together to obtain the proof of Theorem  \ref{T:HKLRP}

\subsection{Estimates on the Growth of the Walk}
Requirement \eqref{I:5} is most conveniently derived using the "the environment seen from the particle".
Let $X^{\omega}_t$ denote the random walk trajectory generated by $P^\omega(x, y)$.  
Recall the shift operation, $\tau_x$ from Section $2$.  By our assumption of translation invariance of the $p_{x, y}$, $\mu$ is clearly translation invariant for all the shifts.  The Kolmogorov $0-1$ law implies that $\mu$ is ergodic with respect to the collection of shifts $\{\tau_x\}_{x \in \mathbb Z^d}$.
Given an initial environment $\omega$, $ \tau_{X^{\omega}_t}: \Omega \rightarrow \Omega$ defines a stochastic map on the space of environments.  Let $\omega_{t}:= \tau_{X^{\omega}_t}(\omega)$, with initial environment $\omega_0=\omega$. It is clear that $\omega_t$ is Markov, since the underlying random walk is.  

Further, given an environment $\omega$, let $d^\omega(0)$ denote the degree of $\omega$ at $0$.
Let $\textd \BbbP(\omega) = \frac{d^\omega(0)}{\E_{\mu}[d^\omega(0)]} \textd \mu(\omega)$ and let us introduce the Hilbert space $L^2(\BbbP)=\{ f: \Omega \rightarrow \mathbb R: \E_{\BbbP}(f^2) < \infty\}$, with inner product $\langle f, g\rangle :=\int \textrm{d} \BbbP(\omega) f(\omega) g(\omega)$.  Note here that since $s \in (d, \infty)$, $\E_{\mu}[d^{\omega}(0)] < \infty$.  

It follows that the operator $A_t f(\omega):= f(\omega_t)$ is self adjoint since the underlying walk is reversible with respect to the un-normalized measure $d^{\omega}(x)$.  Let $Q(\omega, \textd \omega')$ denote the transition kernel for $\omega_t$ going from $\omega$ to $\omega'$ and let $D_t = \max_{0\leq u \leq t} |X_u|$ denote the diameter of the walk at time $t$.
\begin{lemma}
\label{L:growth}
Let $1/(s-d) <p$.  Then for either the discrete or continuous time process, there exists a constant $c$ so that for any $x \in \Z^d$, 
\be
P^{\omega}_{x}\left(D_t > c t^{p+1} \text{ infinitely often }\right) = 0
\ee
$\mu \: a.s.$
Moreover, there exist constants $ c_1, c_2, c_3>0$ such that for any $T, \lambda, p,r>0$ with $p$ as above and $r< s-d$,
\be
 \BbbP(\{P^{\omega}_0(\exists t \leq T: |D_t| \geq c_1T^{p+1}) > c_2/T^{\lambda} \}) \leq c_3T^{\lambda+1-pr}.
\ee
\end{lemma}
\begin{proof}
Discrete Time Case:
Let us consider the increment of the walk at time $n$:
\be
I_n = |X_n-X_{n-1}|.
\ee
Obviously 
\be
D_n \leq \sum_{j=1}^{n} I_j
\ee
Now $I_1 \in L^{r}(\BbbP)$ for any $r < s-d$.
Let us begin by noting that by stationarity of the "environment viewed from the particle" process
\be
\BbbP(\exists t \leq n: I_t \geq n^p) \leq n^{1-pr}.
\ee 
But
\be
\BbbP(\exists t \leq n: I_t \geq n^p) = \E_{\BbbP}(P^{\omega}_0(\exists t \leq n: I_t \geq n^p) )
\ee
so we conclude that
\be
\sum_{k=1}^{\infty} P^{\omega}_0(\exists t \leq 2^k: I_t \geq 2^{kp})
\ee
is summable $\BbbP$, and hence $\mu \: a.s.$ as long as $pr>1$.

Moreover, if we require quantitative bounds, these can be achieved
using Markov's inequality:
\be
\BbbP( \{P^{\omega}_0(\exists t \leq n: I_t \geq n^p) > 1/n^{\lambda} \}) \leq n^{\lambda+1-pr}
\ee

The continuous time case now follows easily by a coupling argument.
\end{proof}

\subsection{The Upper Bound}
The upper bound in Theorem \ref{T:HKLRP} now follows by combining Corollary \ref{C:Gap}, Lemma \ref{L:growth} and Lemma \ref{L:AHK}.

For the sake of completeness we will be explicit.  Consider the event $\{x \in \CC^{\infty}(\omega)\}$.
Hypothesis \eqref{I:6} of Lemma \ref{L:AHK} sets the scale for our calculations.  It is a general fact that, for \textit{any} infinite connected graph $\GG$, $P_t(z, z) \leq c/\sqrt{t}$ \cite{CG} for some constant $c$   independent of $z \in \GG$.
Thus if $V_t$ satisfies \eqref{I:6}, then $V_t \geq c' \sqrt{t}/\log t$.

Let $\kappa, \eta>0$ be fixed.
By Lemma \ref{L:growth}, we can find $p, r$ and a random variable $T_x(\omega)$ so that
\be
\label{g}
P^{\omega}_x(\exists u \leq t: \|D_t\|_2 \geq c_1t^{p}) < c_2/t^{d/(s-d)+ \kappa}
\ee
for all $t\geq \tilde T_{x}(\omega)$ and moreover
\be
\mu(\tilde T_{x}(\omega)>k)< C(\eta, \kappa) k^{-\eta}.
\ee
Choose $\epsilon=1/(4p+ 3d)$.  By Corollary \ref{C:Gap} and the translation invariance of $\mu$,  we can find a positive random variable $T^*_x(\omega)>0$ so that the event
\be
\tau_x \left(B(\epsilon, \lceil 2 c_1t^{p} \rceil)\right)
\ee 
holds for all $t \geq T^*_{x, \epsilon}$ with
\be
\label{sg}
\mu(T^*_{x, \epsilon}>k)< C(\epsilon) e^{-c(\epsilon) \log^2 k}.
\ee
Finally, by Lemma \ref{L:Deg} choosing $\delta =\eta+1$, we can find a random variable $T^{**}_x(\omega)$ so that for all $t \geq T^{**}_x(\omega)$, $\textrm{Deg}(B_x(2c_1t^{p})) \leq 2 c_1 p \log t$ and 
\be
\mu(T^{**}_x(\omega)> k) < C(\eta) k^{-\eta}
\ee

Let $T_x(\omega)= \max(\tilde T_{x}(\omega), T^*_{x, \epsilon}(\omega), T^{**}_x(\omega))$
In the notation of Lemma \ref{L:AHK},  for each $t \geq T_x(\omega)$, let
\begin{align*}
B(t):=&B_x(2 c_1t^{p})\\
\end{align*}
If $t \geq T_x(\omega)$ then by our choices above, we may take $\scrO_t$ to be the minimal cover of 
$B_x(2 c_1t^{p})$ by boxes of side length $\lfloor V^{1/d}_t \rfloor$ and 
\begin{align*}
\PP_t:=&\{\scrG(0, \Lambda): \Lambda \in \scrO_t\}\\
\gamma:=& (s-d)/d
\end{align*}

Thus, for any $t \geq T_x(\omega)$, our choices verify Hypotheses \eqref{I:1}, \eqref{I:2}, \eqref{I:3}, \eqref{I:4}, \eqref{I:5} of Lemma \ref{L:AHK} with $\tilde{\delta}_1= \delta_2$ from Theorem \ref{T:gap} and $\tilde{\delta}_2= 1$.
\qed

\bibliographystyle{plain}
\bibliography{Bib}

\end{document}